\documentclass[10pt]{article}
\usepackage{makeidx}
\usepackage{multirow}
\usepackage{multicol}
\usepackage{graphicx}
\usepackage{epstopdf}
\usepackage{ulem}
\usepackage{hyperref}
\usepackage{amsmath}
\usepackage{amssymb}
\usepackage{setspace}
\usepackage{color}
\usepackage{float}
\author{Aishwarjya Gogoi}
\title{}
\usepackage[paperwidth=612pt,paperheight=792pt,top=72pt,right=72pt,bottom=72pt,left=72pt]{geometry}

\makeatletter
	\newenvironment{indentation}[3]%
	{\par\setlength{\parindent}{#3}
	\setlength{\leftmargin}{#1}       \setlength{\rightmargin}{#2}%
	\advance\linewidth -\leftmargin       \advance\linewidth -\rightmargin%
	\advance\@totalleftmargin\leftmargin  \@setpar{{\@@par}}%
	\parshape 1\@totalleftmargin \linewidth\ignorespaces}{\par}%
\makeatother 

\begin{document}
\begin{doublespace}
\begin{center}
\textbf{{\LARGE Stability evaluation of approximate Riemann solvers using the direct Lyapunov method}} \\
A. Gogoi$^*$, J. C. Mandal$^*$ and A. Saraf$^{**}$\\
$^*$Department of Aerospace Engineering, Indian Institute of Technology, Bombay, Mumbai, 400076, India \\
$^{**}$ Aeronautical Development Agency, Vimanapura, Bengaluru, India 560093.
\end{center}
\begin{center}
	\uppercase{ABSTRACT}
\end{center}

The paper presents a new approach of stability evaluation of the approximate Riemann solvers based on the direct Lyapunov method. The present methodology offers a detailed understanding of the origins of numerical shock instability  in the approximate Riemann solvers. The pressure perturbation feeding the density and transverse momentum perturbations is identified as the cause of the numerical shock instabilities in the complete approximate Riemann solvers, while the magnitude of the numerical shock instabilities are found to be proportional to the magnitude of the pressure perturbations. A shock-stable HLLEM scheme is proposed based on the insights obtained from this analysis about the origins of numerical shock instability in the approximate Riemann solvers.  A set of numerical test cases are solved to show that the proposed scheme is free from numerical shock instability problems of the original HLLEM scheme at high Mach numbers.


\section{Introduction }
The Godunov-type approximate Riemann solvers \cite{godu, toro-book-99, toro-book-09} are a popular method of computing the convective fluxes in the Euler and Navier-Stokes equations. The approximate Riemann solvers can be broadly classified into incomplete and complete solvers. In the incomplete solvers, a simpler wave structure that does not include all the waves in the Riemann problem is used in evaluating the numerical flux while all the non-linear and linearly degenerate waves present in the Riemann problem are included in the complete solvers. Examples of incomplete approximate Riemann solvers are the HLL \cite{hll} and HLLE \cite{einf1} schemes while the Roe \cite{roe}, Osher \cite{osher}, HLLEM \cite{einf2} and HLLC \cite{toro-hllc} schemes are examples of complete approximate Riemann solvers. The HLL-CPS scheme proposed by Mandal and Panwar \cite{mandal}, the HLLCM and HLLS schemes proposed by Shen et al. \cite{shen}, and the HLLEC and HLLES schemes proposed by Xie et al. \cite{xie-hllem} can be categorized as incomplete solvers as these schemes exclude one of the linearly degenerate intermediate waves present in the Riemann problem. The Roe scheme \cite{roe} is based on the Riemann solution of a special linearized version of Euler equations. The Roe scheme is capable of resolving the isolated shock and linearly degenerated discontinuities exactly. However, the original Roe scheme is not positivity preserving and tends to produce spurious expansion shocks under certain conditions. Various types of entropy fixes \cite{harten} have been proposed to overcome the spurious expansion shock problem of the Roe scheme. The HLL scheme proposed by Harten, Lax, and van Leer \cite{hll} assumes a two-wave structure of the Riemann problem comprising the fastest left and right-running acoustic waves. The HLL scheme with the wave-speed estimates of Einfeldt \cite{einf1} satisfies positivity and entropy conditions desirable for a numerical scheme. Due to its desirable positivity-preserving and entropy-satisfying properties, the HLL framework has been utilized to develop complete approximate Riemann solvers like HLLEM \cite{einf2} and HLLC \cite{toro-hllc}. In the HLLEM scheme proposed by Einfeldt et al. \cite{einf2}, anti-diffusion terms are added to the HLLE scheme corresponding to the linearly degenerate waves that are not resolved in the HLL scheme. The accuracy of the HLLEM scheme is similar to the Roe scheme and the scheme retains the positivity-preserving and entropy-satisfying properties of the HLLE scheme. In the HLLC scheme proposed by Toro et al. \cite{toro-hllc}, the missing contact and shear waves were restored into the HLL scheme by utilizing the Rankine-Hugoniot jump conditions. The HLLC scheme with the wave speed estimates of Batten et al. \cite{batt} has emerged as perhaps the most preferred Riemann solver due to its ability to resolve the linearly degenerate waves. 

It may be noted that for high-speed flows, the complete approximate Riemann solvers like the Roe \cite{roe}, Osher \cite{osher}, HLLEM \cite{einf2}, and HLLC \cite{toro-hllc} schemes are susceptible to numerical shock instabilities like odd-even decoupling, kinked Mach stem, carbuncle phenomenon, and low-frequency post-shock oscillations in slowly moving shock. On the other hand, the HLLE scheme \cite{hll, einf1}, which does not resolve the linearly degenerate intermediate waves, is free from numerical shock instabilities. The cause of numerical shock instabilities in the complete approximate Riemann solvers has been widely investigated ever since they were first presented by Peery and Imlay \cite{peery} in their hypersonic blunt body computations using the Roe scheme. Quirk \cite{quirk} has carried out linear perturbation analysis and has found that the schemes in which the pressure perturbations feed the density perturbations will be afflicted by the odd-even decoupling problem. Sanders et al. \cite{sand} have observed that the instability is a result of inadequate cross-flow dissipation offered by the strictly upwind schemes and proposed a multi-dimensional, upwind dissipation to eliminate instability in the Roe scheme. Pandolfi and D'Ambrosio \cite{pand} have carried out a more refined perturbation analysis by including the velocity perturbations and have observed that the schemes in which the density perturbations are fed by the pressure perturbations and the shear velocity perturbations are undamped shall suffer from severe carbuncle phenomenon. Pandolfi and D'Ambrosio \cite{pand} have also observed that the scheme in which the density perturbations are damped while the shear velocity perturbations are not damped shall suffer from mild carbuncle phenomenon. Liou \cite{liou1} has conjectured that the schemes that suffer from shock instability have a pressure jump term in the numerical mass flux while those free from shock instability do not contain a pressure jump term in the mass flux. Based on the conjecture of Liou \cite{liou1}, Kim et al. \cite{kim} have proposed a shock-stable Roe scheme by balancing the feeding and damping of the density and pressure terms in the mass flux. However, the conjecture of Liou is not universally accepted and the AUSM family of schemes \cite{ausm1} are often cited as counterexamples to this conjecture as these schemes do not contain pressure jump terms in the mass flux but exhibit carbuncle phenomenon under conditions like highly stretched grids. Gressier and Moschetta \cite{gress1} have carried out linear stability analysis of the perturbation evolution equations \cite{quirk, pand} and have shown that any upwind scheme that resolves the contact discontinuity exactly is marginally stable. Gressier and Moschetta \cite{gress1} have also proposed the conjecture that strict stability and exact resolution of contact discontinuities are not compatible. Dumbser et al. \cite{dumb-matrix} have developed a useful matrix stability method to evaluate the stability of the numerical schemes and observed that the origin of the carbuncle phenomenon is localized in the upstream region of the shock and the location of the intermediate point in the numerical shock structure plays a very key role in the onset of instability. If the intermediate state is sufficiently close to the downstream state, the carbuncle phenomenon will not appear even at very high upstream Mach numbers. The closer the intermediate state is to the upstream state, the greater the susceptibility of the scheme to encounter the carbuncle phenomenon. The matrix stability method of Dumbser et al. is perhaps the most popular method of evaluating the stability of upwind schemes. The effect of numerical shock structure on shock instability have been further investigated by Gogoi et al. \cite{gogoi} for the HLL-family schemes. The exact resolution of the shear wave has been  identified by Shen et al. \cite{shen} as the cause of the numerical instabilities in the HLLC scheme. Shen et al. \cite{shen} have shown that the HLLC scheme minus the shear wave (named HLLCM) is free from numerical shock instability while the HLL scheme with the shear wave (named HLLS) is susceptible to numerical shock instability problems. Shen et al. \cite{shen} have proposed a hybrid HLLC-HLLCM scheme that switches over to the HLLCM scheme in the vicinity of the shock while the HLLC scheme is used elsewhere. Similarly, Xie et al. \cite{xie-hllem} have identified the exact resolution of the shear wave as the cause of numerical shock instabilities in the HLLEM Riemann solver. Xie et al. \cite{xie-hllem} have shown that the HLLEC scheme, which is the HLLEM scheme with anti-diffusion for the contact wave, is free from numerical shock instability problems while the HLLES scheme, which is the HLLEM scheme with anti-diffusion for the shear wave, is susceptible to numerical shock instability problems. Xie et al. \cite{xie-hllem} have proposed an HLLEM scheme that reduces the contribution of the shear wave in the vicinity of the shock to overcome the numerical shock instability problems in the original HLLEM scheme. However, Kemm \cite{kemm} has shown that the resolution of the entropy or the contact wave in the HLLEM Riemann solver also contributed to the numerical shock instabilities, although its contribution to the numerical shock instabilities is less as compared to the shear wave. Simon and Mandal \cite{san-hllem} have shown that the HLLEMS scheme of Xie et al. \cite{xie-hllem} is not free from numerical shock instability problems. Simon and Mandal \cite{san-hllem} have shown that suppression of the anti-diffusion terms of both the contact and shear waves is more effective in curing the numerical shock instability problems in the HLLEM Riemann solver. Fleischmann et al. \cite{fleis1, fleis2} have identified the inappropriate numerical dissipation in the transverse direction of a grid-aligned shock to be the cause of the numerical shock instability in the Godunov-type approximate Riemann solvers like the Roe scheme and the HLLC scheme. The inappropriate numerical dissipation in the transverse direction is attributed to the low Mach number in the transverse direction of the grid-aligned shock. Fleischmann et al. \cite{fleis1} have proposed a reduction of the overall numerical dissipation in the classical Roe scheme and the component-wise local Lax-Fredrichs (cLLF) scheme to suppress the shock instability and the reduction in the overall numerical dissipation is achieved by modifying the non-linear acoustic eigenvalues. Fleischmann et al. \cite{fleis2} have observed that similar modification to the non-linear eigenvalues in the HLLC scheme leads to a dissipative upwind scheme. Fleischmann et al. \cite{fleis2} have proposed a centred HLLC scheme comprising central-difference terms and numerical dissipation terms, and proposed modification of the non-linear eigenvalues in the numerical dissipation terms to suppress shock instability. Fleischmann et al. \cite{fleis2} have also demonstrated reduced overall numerical dissipation with the modification of the non-linear eigenvalues and have shown accurate results over a cylinder at low Mach numbers.

It is observed that the investigation into the origins of numerical shock instabilities in approximate Riemann solvers is ongoing. Moreover, the existing methods for assessing the stability of numerical schemes rely on linearising the fluxes, and their outcomes do not consistently align with numerical results. For instance, linear perturbation analysis suggests stability for the Roe scheme with Harten's entropy fix \cite{harten} applied to linearly degenerate eigenvalues across all entropy fix values. However, numerical findings indicate that stability of the Roe scheme is achieved only beyond a specific threshold value of entropy fix. Consequently, there is a perceived need for more sophisticated stability assessment techniques to accurately pinpoint the causes of numerical shock instabilities in approximate Riemann solvers.
In this paper, we evaluate the stability of the HLL-family schemes using the direct Lyapunov method. Also a shock-stable HLLEM scheme is proposed based on the stability evaluation. The paper is organized into eight sections. In section 2, the governing equations are described. In section 3, finite volume discretization is described. In section 4, linearized perturbation analysis of the HLL-family schemes is presented while in section 5, stability analysis of the HLL-family scheme is performed using the direct Lyapunov method. In section 6, the formulation of a shock-stable HLLEM scheme is presented and the numerical results of the shock-stable HLLEM scheme are presented in section 7. Finally, conclusions are drawn in section 8.

\section{Governing equations}
The Euler equations for the two-dimensional inviscid compressible flow in differential, conservative form is
\begin{equation} \label{ee-dc}
\dfrac{\partial{}\mathbf{\acute{U}}}{\partial{}t}+\dfrac{\partial{}\mathbf{\acute{F}(\acute{U})}}{\partial{}x}+\dfrac{\partial{}\mathbf{\acute{G}(\acute{U})}}{\partial{}y}=0
\end{equation}
where $\mathbf{\acute{U}}$, $\mathbf{\acute{F}(\acute{U})}$, and $\mathbf{\acute{G}(\acute{U})}$ are the vector of conserved variables, x-directional and y-directional fluxes respectively, and can be written as
\begin{equation}
\mathbf{\acute{U}}=\left[\begin{array}{c} \rho \\ \rho{}u \\ \rho{}v \\ \rho{}E \end{array}\right] \hspace{1cm} \mathbf{\acute{F}(\acute{U})}=\left[\begin{array}{c} \rho{}u \\ \rho{}u^2+p \\ \rho{}uv \\ (\rho{}E+p)u \end{array}\right]  \hspace{1cm} \mathbf{\acute{G}(\acute{U})}=\left[\begin{array}{c} \rho{}v \\ \rho{}uv \\ \rho{}v^2+p \\ (\rho{}E+p)v \end{array}\right] 
\end{equation}
where $\rho$, $u$, $v$, $p$, and $E$ stand for density, x-directional and y-directional velocities in global coordinates, pressure, and specific total energy. The system of equations can be closed by the equation of state
\begin{equation}
p=(\gamma-1)\left(\rho{}E-\frac{1}{2}\rho(u^2+v^2)\right)
\end{equation} 
where $\gamma$ is the ratio of specific heat. In this paper, a calorically perfect gas with $\gamma=1.4$ is considered. The speed of sound is $ a=\sqrt{\dfrac{\gamma{}p}{\rho}}$. The corresponding integral form  of the governing equations is
\begin{equation}\label{ee-ic}
\frac{d}{dt}\int_{\Omega{}}\mathbf{\acute{U}}d\Omega{}+\oint_{\partial{}\Omega{}}\mathbf{(\acute{F},\acute{G}).n}dS=0
\end{equation}
where $\partial{}\Omega$ is the boundary of the control volume $\Omega$ (area in case of two dimensions) and  $\mathbf{n}$ denotes the outward pointing unit normal vector to the surface $\partial\Omega$, $dS$ is the length of  the edge.

\section{Finite volume discretization} The finite volume discretization of equation (4) for a structured, quadrilateral mesh can be written as
\begin{equation} 
\dfrac{d\mathbf{\acute{U}}_i}{dt}=-\dfrac{1}{|\Omega|_i}\sum_{l=1}^4[\mathbf{(\acute{F},\acute{G})}_l.\mathbf{n}_l]\Delta{}s_l
\end{equation}
where $\mathbf{\acute{U}}_i$ is the cell averaged conserved state vector, $\mathbf{n}_l$ denotes the unit normal vector and $\Delta{}s_l$ denotes the length of each interface. The rotational invariance property of the Euler equations is utilized to express the flux $\mathbf{(\acute{F},\acute{G})}_l.\mathbf{n}_l$ as 
\begin{equation} \label{fv}
\dfrac{d\mathbf{\acute{U}}_i}{dt}=-\dfrac{1}{|\Omega|_i}\sum_{l=1}^{4}(T_{il})^{-1}\mathbf{F}(T_{il}\mathbf{U}_i,T_{il}\mathbf{U}_l)\Delta{}s_l
\end{equation}
where $\mathbf{F}(T_{il}\mathbf{U}_i,T_{il}\mathbf{U}_l)$ is the inviscid face normal flux vector, $T_{il}$ is the rotation matrix and ${T_{il}}^{-1}$ is its inverse at the edge between cell $i$ and its neighbour $l$. 
\begin{equation}
T_{il}=\left(\begin{array}{ccccc}1& 0& 0& 0\\ 0& n_x& n_y& 0 \\  0& -n_y& n_x& 0 \\ 0& 0& 0& 1\end{array}\right) \hspace{1cm} \text{and} \hspace{1cm} T_{il}^{-1}=\left(\begin{array}{ccccc}1& 0& 0& 0\\ 0& n_x&- n_y& 0 \\  0& n_y& n_x& 0 \\ 0& 0& 0& 1\end{array}\right)
\end{equation}
where $\mathbf{n}=(n_x,n_y)$ is the outward unit normal vector on the edge between cells $i$ and $l$. 

In the next section, the HLL-family of approximate Riemann solvers for evaluation of the face normal flux $\mathbf{F}(T_{il}\mathbf{U}_i,T_{il}\mathbf{U}_l)$ shown in equation (\ref{fv}) is described. 
\section{HLL family flux function }
\subsection{The HLL approximate Riemann solver}
The HLL scheme proposed by Harten, Lax, and van Leer \cite{hll} assumes a simplified wave structure of the Riemann problem comprising the fastest left running and right running waves acoustic waves. The approximate Riemann solution in the region bounded by the two acoustic waves is 
\begin{equation}  \label{hll-state-variable}
\mathbf{U}_{*,HLL}=\dfrac{S_R\mathbf{U}_R-S_L\mathbf{U}_L+\mathbf{F}_L-\mathbf{F}_R}{S_R-S_L}
\end{equation}
where $S_R$ and $S_L$ are the fastest right and left running wave speeds, $\mathbf{U}_R=\mathbf{T}_{il}\mathbf{U}_i$ and $\mathbf{U}_L=\mathbf{T}_{il}\mathbf{U}_l$ are the face normal state variables to the right and left of the interface, $\mathbf{F}_R=\mathbf{F}(\mathbf{U}_R)$ and $\mathbf{F}_L=\mathbf{F}(\mathbf{U}_L)$ are the face normal fluxes to the right and left of the interface.
Based on the expression for the state variables, the HLL flux at the interface of two cells is given by \cite{hll}
\begin{equation} \label{flux-hll}
\mathbf{F}(\mathbf{U}_L,\mathbf{U}_R)=\left \lbrace \begin{array}{cl} \mathbf{F}_L & \text{if} \hspace{3mm}0\le S_L\\ 
\dfrac{S_R\mathbf{F}_L-S_L\mathbf{F}_R+S_RS_L(\mathbf{U}_R-\mathbf{U}_L)}{S_R-S_L} & \text{if} \hspace{3mm} S_L\le 0 \le S_R \\  \mathbf{F}_R & \text{if} \hspace{3mm} 0\ge S_R \end{array}\right.
\end{equation}
The wave speed estimates of Einfeldt \cite{einf1} are 
\begin{equation} \label{wavespeed-einf}
S_L=min(u_{nL}-a_L,\tilde{u}_n-\tilde{a}) \hspace{1cm} S_R=max(u_{nR}+a_R,\tilde{u}_n+\tilde{a}) 
\end{equation}
where $u_{nL}$ and $u_{nR}$ are the normal velocities to the left and right of the interface, $a_L$ and $a_R$ are the speed of sound to the left and right of the interface, $\tilde{u}_n$ and $\tilde{a}$ are the Roe-averaged normal velocity and speed of sound at the cell interface. The velocity normal to the interface is $u_n=un_x+vn_y$. The HLL scheme with the wave speed estimates of Einfeldt (\ref{wavespeed-einf}) is commonly referred to as the HLLE scheme. The HLL/HLLE flux can also be written as
\begin{equation} \label{flux-hlle}
\mathbf{F}(\mathbf{U}_R,\mathbf{U}_L)=\dfrac{S_R\mathbf{F}_L-S_L\mathbf{F}_R}{S_R-S_L}+\dfrac{S_RS_L}{S_R-S_L}(\Delta{}\mathbf{U})
\end{equation} 

\begin{indentation}{80pt}{0pt}{0pt} where ${\Delta{}\mathbf{U}=\mathbf{U}_R-\mathbf{U}_L}$ and \end{indentation}
\begin{equation} \label{wavespeed-hll}  
S_L=min(0,\;u_{nL}-a_L,\;\tilde{u}_n-\tilde{a}) \hspace{1cm} S_R=max(0,\;u_{nR}+a_R,\;\tilde{u}_n+\tilde{a})
\end{equation}

\subsection {The HLLEM approximate Riemann solver}
The HLLEM flux is the HLL/HLLE flux with additional anti-diffusion terms for the linearly degenerate intermediate waves that are not included in the HLL/HLLE scheme. The HLLEM flux is given by \cite{einf2}
\begin{equation} \label{hll-type-flux}
\mathbf{F}(\mathbf{U}_R,\mathbf{U}_L)=\dfrac{S_R\mathbf{F}_L-S_L\mathbf{F}_R}{S_R-S_L}+\dfrac{S_RS_L}{S_R-S_L}(\Delta{}\mathbf{U}-B\Delta{}\mathbf{U})
\end{equation} 
where $B\Delta{}\mathbf{U}$ is the anti-diffusion term for the linearly degenerate intermediate waves. The wave speeds are the same as that of the HLLE scheme shown in equation (\ref{wavespeed-hll}). Einfeldt et al. \cite{einf2} have shown that with the above choice of wave speeds, the HLLEM scheme is entropy-satisfying and positively conservative \cite{einf2} for a perfect gas. The anti-diffusion term $B\Delta{}\mathbf{U}$ in the HLLEM scheme comprises anti-diffusion coefficients, wave strength, and right eigenvectors corresponding to the linearly degenerate characteristic fields of the flux Jacobian matrix. The anti-diffusion term is \cite{einf2,kim}
\begin{equation} \label{bdq-roe-hllem}
B\Delta{}\mathbf{U}=\delta_2\tilde{\alpha_2}\tilde{\mathbf{R}}_2+ \delta_3\tilde{\alpha_3}\tilde{\mathbf{R}}_3=\delta_2\left(\Delta{}\rho-\dfrac{\Delta{}p}{\tilde{a}^2}\right) \left[1,\;\tilde{u}_n,\;\tilde{u}_t,\;\frac{1}{2}(\tilde{u}_n^2+\tilde{u}_t^2)\right]^T+\delta_3\tilde{\rho}\Delta{}u_t\left[0,\;0,\;1,\;\tilde{u}_t\right]^T
\end{equation}
where $\Delta(.)=(.)_R-(.)_L$, $\delta_2$ and  $\delta_3$ are the anti-diffusion coefficients for the contact and shear waves which are not resolved in the classical HLL scheme. The anti-diffusion coefficients are defined as $\delta{}_2=\delta{}_3=\frac{\tilde{a}}{\tilde{a}+|\bar{u}|}$ where $|\bar{u}|$ is the approximate speed of the contact discontinuity given by $|\bar{u}|=\frac{S_R+S_L}{2}$ \cite{einf2,kim}. Park and Kwon \cite{park} observed that the HLLEM scheme of Einfeldt \cite{einf2} is not as accurate as the Roe scheme in the presence of a contact discontinuity and have proposed the Roe-averaged normal velocity as the speed of the contact discontinuity to improve accuracy of the HLLEM scheme. Thus, for the HLLEM scheme with the modification of Park and Kwon, $\delta{}_2=\delta{}_3=\frac{\tilde{a}}{\tilde{a}+|\tilde{u}_n|}$. This expression of the anti-diffusion coefficients shall be used in the present work for the HLLEM scheme.

\section{Linear perturbation analysis}
Linear perturbation analysis of the numerical scheme is useful for determining the stability of the numerical scheme for high-speed flow problems. Linear perturbation analysis of various HLL-family schemes are carried out like Quirk \cite{quirk} and Pandolfi-D'Ambrosio \cite{pand} for the case of a grid aligned-shock with zero normal velocity across the interface. The evolution of the density, shear velocity, and pressure perturbations in the HLL-family scheme is obtained from the linear perturbation analysis. A saw-tooth type perturbation is considered and the normalized flow properties are described as follows
\begin{equation}
\rho=1\pm \hat{\rho}, \hspace{1cm} u=u_0\pm \hat{u}, \hspace{1 cm} v=0, \hspace{1cm}  p=1\pm\hat{p}
\end{equation}
where $\hat{\rho}$, $\hat{u}$ and $\hat{p}$ are the density, shear velocity, and pressure perturbations. It may be noted that here $u$ is the shear velocity and $v$ is the normal velocity. The evolution of the density, shear velocity, and pressure perturbations in the various HLL-family schemes are shown in Table \ref{table-lpa}. The HLL-family schemes analysed  here are the original HLLE scheme, the complete or all-wave Roe, HLLEM and HLLC schemes, the contact resolving but shear dissipative HLL-CPS, HLLCM and HLLEC schemes, and the shear resolving but contact dissipative HLLS and HLLES schemes. It has been mentioned by Quirk \cite{quirk} that the schemes in which the density perturbations are fed by the pressure perturbations are susceptible to numerical shock instabilities. The all-wave Roe, HLLEM and HLLC schemes are known to be carbuncle-prone schemes, and the pressure perturbations feed the density perturbations in these schemes as seen in Table \ref{table-lpa}. The numerical shock instabilities in these all-wave schemes are consistent with the conjecture of Quirk. It can be seen from the above table that the pressure perturbations feed the density perturbation in the HLL-CPS, HLLCM and HLLEC schemes just like the Roe, HLLEM, and HLLC schemes \cite{quirk, pand}. However, the shear velocity perturbations are damped in the HLL-CPS, HLLCM and HLLEC, like the HLLE scheme. These contact-capturing schemes have been claimed by their authors to be carbuncle-free schemes, but are expected to be carbuncle-prone according to the conjecture of Quirk \cite{quirk} and Liou\cite{liou1}. In recent years, carbuncle has been reported in these schemes by Simon and Mandal \cite{san-hllem} and Gogoi et al.\cite{gogoi}. In the shear-resolving HLLS and HLLES schemes, the density and pressure perturbations are damped, while the shear velocity perturbations are not damped. The HLLS and HLLES schemes have been shown by Shen et al. \cite{shen} and Xie et al. \cite{xie-hllem} to be prone to numerical shock instabilities and the exact resolution of the shear wave has been attributed as the cause of numerical shock instabilities in these schemes. The conjectures of Quirk \cite{quirk} and Liou \cite{liou1} are unable to explain the numerical shock instabilities in these shear wave resolving schemes whereas the instability in these schemes are consistent with the observation of Pandolfi and D'Ambrosio \cite{pand}.

\begin{table}[H]
	\caption{Result of linear perturbation analysis of HLL-family schemes}
	\label{table-lpa}
\begin{center}

\begin{doublespace}
\begin{tabular}{|c|c|c|c|c|}
\hline Serial No & Scheme & $\hat{\rho}^{n+1}=$ & $\hat{u}^{n+1}=$  & $\hat{p}^{n+1}=$\\
\hline 1 & HLLE & $\hat{\rho}^n(1-2\nu)$ & $\hat{u}^n(1-2\nu)$  & $\hat{p}^n(1-2\nu)$\\
\hline 2 &Roe/ HLLEM/HLLC & $\hat{\rho}^n-\dfrac{2\nu}{\gamma}\hat{p}^n$ & $\hat{u}^n$ & $\hat{p}^n(1-2\nu)$\\
\hline 3 & HLL-CPS & $\hat{\rho}^n-\dfrac{2\nu}{\gamma}\hat{p}^n$ & $\hat{u}^n(1-\dfrac{2\nu}{\gamma})$& $\hat{p}^n(1-2\nu)$ \\
\hline 4 & HLLCM and HLLEC & $\hat{\rho}^n-\dfrac{2\nu}{\gamma}\hat{p}^n$ & $\hat{u}^n(1-2\nu)$& $\hat{p}^n(1-2\nu)$ \\
\hline 5 & HLLS and HLLES & $\hat{\rho}^n(1-2\nu)$ & $\hat{u}^n$& $\hat{p}^n(1-2\nu)$ \\
\hline
\end{tabular}
\end{doublespace}
\end{center}

\begin{indentation}{0pt}{0pt}{60pt}

where 

$\nu$ is the linearized Courant number, 

$\gamma$ is the ratio of specific heat, 

$n$ and $n+1$ are the time steps.
\end{indentation}

\end{table}

It becomes evident that the conjecture proposed by Quirk \cite{quirk}, utilizing the linear perturbation analysis method, falls short in establishing the stability of some of certain schemes within the HLL-family. Consequently, there arises a need for a more sophisticated analysis to assess the stability of the approximate Riemann solvers. In the subsequent section, we present a more refined analysis based on the direct Lyapunov method focussing on the HLL-family of approximate Riemann solvers.

\section {Refined analysis of HLL-type schemes}
The Lyapunov method stands out as a widely embraced approach for evaluating the stability of non-linear system of differential equations. It is important to highlight that the continuous Euler equations, being inherently non-linear, and the HLL-family of flux functions, also characterized by non-linearity, underscores the necessity  of employing non-linear analysis. Thus a deliberate effort is made in this study to apply non-linear analysis to assess the stability of the HLL-family schemes.
\subsection{Reduced or first Lyapunov method for discrete Euler equations}
The reduced or first Lyapunov method is a commonly used approach to assess the stability of non-linear systems, often preceding the use of the direct Lyapunov method. This method is specifically applied to the linearized form of the original non-linear system. In the context of the discrete Euler equations, the reduced Lyapunov method corresponds to the linearized perturbation method proposed by Quirk. Notably, the stability analysis conducted by Gressier and Moschetta \cite{gress1} for the perturbation evolution equations bears similarity to the reduced Lyapunov method. The linearized perturbation system proposed by Quirk for the discrete Euler equations can be expressed as
\begin{equation} U^{n+1}=f(U^n) \end{equation}
where  $U=[\rho, u, v, p]^T, f$ is continuously differentiable in the neighbourhood of the origin, and $n, n+1$ are the successive time steps. Let $J=\left[\dfrac{\partial{}f}{\partial{}U^n}\right]_{U^n=U_0}$ be the Jacobian of the system evaluated at this equilibrium. If all the eigenvalues of $J$ are strictly less than unity in the absolute value, then the system is asymptotically stable about its zero equilibrium. If one or more eigenvalues of the system Jacobian $J$ is more than unity in absolute value, the system is unstable about its zero equilibrium. If one or more eigenvalues of the system Jacobian is equal to unity in absolute value, the reduced Lyapunov criterion cannot be used to evaluate the stability of the system.

The perturbation equations of the HLLE scheme are 
\begin{equation}
\left(\begin{array}{c}  \rho \\ u \\ p \end{array}\right)^{n+1}=\left[ \begin{array}{cccc} 1-2\nu & 0 & 0 \\  0 &1-2\nu & 0 \\ 0 & 0 & 1-2\nu\end{array}\right]\left(\begin{array}{c}  \rho \\ u \\ p \end{array}\right)^{n}
\end{equation}
where $\nu$ is the linearized CFL number. The eigenvalues or the amplification factors for the HLLE scheme are $1-2\nu, 1-2\nu$, and $1-2\nu$. All the eigenvalues for the HLLE scheme are less than unity and hence, the HLLE scheme is asymptotically stable. 

The perturbation equations of the Roe, HLLEM and HLLC schemes are 
\begin{equation}
\left(\begin{array}{c}  \rho \\ u \\ p \end{array}\right)^{n+1}=\left[ \begin{array}{cccc} 1 & 0 & \dfrac{-2\nu}{a_0^2} \\  0 &1 & 0 \\ 0 & 0 & 1-2\nu\end{array}\right]\left(\begin{array}{c}  \rho \\ u \\ p \end{array}\right)^{n}
\end{equation}
The eigenvalues of the Roe, HLLEM and HLLC schemes are 1, 1, and 1$-2\nu$. It can be seen that two eigenvalues are equal to unity. Hence, the stability of the Roe/HLLEM/HLLC schemes cannot be determined by the reduced Lyapunov method. The Roe, HLLEM, and HLLC schemes have been described by Gressier and Moschetta \cite{gress1} as marginally stable schemes.

The perturbation equations of the HLLCM scheme of Shen et al. \cite{shen} and HLLEC scheme of Xie et al. \cite{xie-hllem} are 
\begin{equation}
\left(\begin{array}{c}  \rho \\ u \\ p \end{array}\right)^{n+1}=\left[ \begin{array}{cccc} 1 & 0 & \dfrac{-2\nu}{a_0^2} \\  0 &1-2\nu & 0 \\ 0 & 0 & 1-2\nu\end{array}\right]\left(\begin{array}{c}  \rho \\ u \\ p \end{array}\right)^{n}
\end{equation}
The eigenvalues of the HLLEC scheme are $1, 1-2\nu,$ and 1$-2\nu$. It can be seen that one eigenvalue is equal to unity. Hence, the stability of the HLLCM and HLLEC schemes cannot be determined by the reduced Lyapunov method.

The perturbation equations of the HLLS scheme of Shen et al. \cite{shen} and HLLES scheme of Xie et al. \cite{xie-hllem} are 
\begin{equation}
\left(\begin{array}{c}  \rho \\ u \\ p \end{array}\right)^{n+1}=\left[ \begin{array}{cccc} 1-2\nu & 0 & 0 \\  0 &1 & 0 \\ 0 & 0 & 1-2\nu\end{array}\right]\left(\begin{array}{c}  \rho \\ u \\ p \end{array}\right)^{n}
\end{equation}
The eigenvalues of the HLLES scheme are $1-2\nu,1,$ and 1$-2\nu$. It can be seen that one eigenvalue is equal to unity. Hence, the stability of the HLLS and HLLES schemes cannot be determined by the reduced Lyapunov method.

\subsection{Direct Lyapunov method} \label{direct-lyapunov}
Due to the inconclusive nature of the linear perturbation method and the reduced Lyapunov method, the direct Lyapunov method will be better in assessing the stability of the HLL-family schemes. The direct Lyapunov method comprises of finding a suitable Lyapunov function $V(\mathbf{U})$ which satisfies the following criterion
\begin{enumerate}
\item $V(\mathbf{U})$ is continuously differentiable.
\item $V(\mathbf{U}) > 0$ in the neighbourhood of the steady state $\mathbf{U}_0$  and
\item $V(\mathbf{U}_0)=0$.
\end{enumerate}
A quadratic Lyapunov function has been proposed by Bastin et al. \cite{bastin} for the continuous, isothermal Euler equations. Like most Lyapunov functions, this Lyapunov function has the unit of energy. In the present work, a similar Lyapunov candidate function is proposed for the discrete Euler equations. In the direct Lyapunov method, the original non-linear equations are considered. A candidate Lyapunov function for the discrete Euler equations with fluxes evaluated by the approximate Riemann solver is taken as
\begin{equation} \label{lf}
V(\mathbf{U})=\dfrac{a_0^2}{\rho{}_0}\hat{\rho}^2+\dfrac{a_0^2}{\rho_0u^2_0}(\hat{\rho{}u})^2+\dfrac{1}{\rho_0a_0^2}\hat{p}^2
\end{equation}
where $\mathbf{U}=[\rho, \rho{}u,\rho{}v, p]^T$ represents the state variables, $\hat{\rho}, \hat{\rho{}u}$, and $\hat{p}$ are the deviation of the states $\rho,\rho{}u$, and $p$ from their steady states $\rho_0, (\rho{}u)_0$ and $p_0$, and $a_0$ is the steady state speed of sound. The deviation of the states can be expressed as
\begin{equation}
\hat{\rho}=\rho(t)-\rho_0, \hspace{1 cm} \hat{\rho{}u}=\rho{}u(t)-\rho_0u_0, \hspace{1 cm} \hat{\rho{}v}=0, \hspace{1cm} \hat{p}=p(t)-p_0
\end{equation}
Here the state variables are taken as $\mathbf{U}=[\rho,\;\rho{}u,\;\rho{}v,\;p]^T$ whereas the state variables in the original Euler equations are $\mathbf{U}=[\rho,\;\rho{}u,\;\rho{}v,\;\rho{}e]^T$. Here, it may be noted that considering pressure as the state variable instead of total internal energy leads to a major simplification of the derivation  without any compromise on the results. The total internal energy comprises pressure and kinetic energy and the terms that contribute to the kinetic energy are already included in the continuity and momentum equations.
\subsubsection{Global asymptotic stability criterion}
For the discrete Euler equations, the global asymptotic stability criterion of Lyapunov can be expressed as
\begin{enumerate}
\item $V(\mathbf{U}_0)=0$.
\item $V(\mathbf{U})^n> 0$  $\forall$ $\mathbf{U}^n \neq 0$ in $\mathcal{D}$.
\item  $\Delta{}V(\mathbf{U})^n:=V(\mathbf{U})^{n+1}-V(\mathbf{U})^n < 0$   $\forall$ $\mathbf{U}^n \neq 0$ in $\mathcal{D}$.
\item $V(\mathbf{U})\to \infty$ as $||\mathbf{U}|| \to \infty$, i.e. $V(\mathbf{U})$ is radially unbounded.
\end{enumerate}
Here the superscript $n$ stands for the time step.

The Lyapunov function shown in equation (\ref{lf}) is positive definite in the neighbourhood of the steady state $\mathbf{U}_0$ and is equal to zero at the steady state. Thus the function can be considered as an appropriate Lyapunov function. The delta change in the Lyapunov function with every time step used for evaluating the stability of the numerical scheme can be written as
\begin{equation} \label {lf1}
\Delta{}V(\mathbf{U}):=V(\mathbf{U})^{n+1}-V(\mathbf{U})^n=\dfrac{a_0^2}{\rho{}_0}(\hat{\rho}^{2,n+1}-\hat{\rho}^{2,n})+\dfrac{a_0^2}{\rho_0u_0^2}(\hat{\rho{}u}^{2,n+1}-\hat{\rho{}u}^{2,n})+\dfrac{1}{\rho_0a_0^2}(\hat{p}^{2,n+1}-\hat{p}^{2,n})
\end{equation}
The results of the direct Lyapunov method for the various HLL-family schemes are presented in the next sections. Here, the analysis approach is similar to the linear perturbation method of Quirk \cite{quirk} and Pandolfi-D'Ambrosio \cite{pand}, with the normal velocity is taken as $v=0$. A saw-tooth-type perturbation in the flow variables is assumed and the normalized flow variables are described as 
\begin{equation}
\rho=\rho_0\pm \hat{\rho}, \hspace{1cm} \rho{}u=(\rho{}u)_0\pm \hat{\rho{}u}, \hspace{1 cm} \rho{}v=0, \hspace{1cm}  p=p_0\pm\hat{p}
\end{equation}
where $\rho_0$, $(\rho{}u)_0$, $(\rho{}v)_0$, $p_0$ are the steady state values of density, transverse momentum, normal momentum, and pressure, $\hat{\rho}$, $\rho{}\hat{u}$, $\rho{}\hat{v}$ and $\hat{p}$ are the density, transverse momentum, normal momentum, and pressure perturbations. It may be noted that here $\rho{}u$ is the transverse momentum velocity and $\rho{}v$ is the normal momentum. Just like the analysis of Quirk \cite{quirk}, it is assumed that the stream-wise flux is balanced, i.e., $(\mathbf{F}_{i+\frac{1}{2},j}-\mathbf{F}_{i-\frac{1}{2},j}=0)$ and only the transverse flux remains. Therefore, the Euler equation (\ref{fv}) reduces to 
\begin{equation}
\mathbf{U}^{n+1}_{i,j}=\mathbf{U}^n_{i,j}-\dfrac{\nu}{a_0+|v_0|}(\mathbf{G}_{i,j+\frac{1}{2}}-\mathbf{G}_{i,j-\frac{1}{2}})
\end{equation} 
where $\mathbf{U}=[\rho,\;\rho{}u,\;\rho{}v,\;\rho{}e]^T$ is the vector of state variables, $\nu=\frac{\Delta{}t(a_0+|v_0|)}{\Delta{}y}$ is the linearized CFL number, superscripts $n$ and $n+1$ are the time steps, and subscripts $i$ and $j$ are the cell indices. The perturbation equations are obtained by cancelling the steady state values of the above equation and can be written as 
\begin{equation}
\mathbf{\hat{U}}^{n+1}_{i,j}=\mathbf{\hat{U}}^n_{i,j}-\dfrac{\nu}{a_0}(\mathbf{G}_{i,j+\frac{1}{2}}-\mathbf{G}_{i,j-\frac{1}{2}})
\end{equation} 
where $\mathbf{\hat{U}}=[\hat{\rho},\;\hat{\rho{}u},\;\hat{\rho{}v},\;\hat{\rho{}e}]^T$ is the vector of the perturbations in the state variables $\mathbf{U}$.
\subsection{Analysis of the HLLE scheme}
The flux difference is obtained from the HLLE flux equations (\ref{fv}, \ref{flux-hlle}) and is  
\begin{equation}
\begin{array}{l}
\mathbf{G}_{i,j+\frac{1}{2}}-\mathbf{G}_{i,j-\frac{1}{2}}= a_0\left[\begin{array}{c} 2\hat{\rho}^n\\  2\hat{\rho{}u}^n \\ \dfrac{2\hat{p}^n}{\gamma-1}+\hat{\rho{}q^2}^n \end{array}\right]
\end{array}
\end{equation}
The continuity perturbation equation can be written as
\begin{equation}
\hat{\rho}^{n+1}=\hat{\rho}^n-2\nu\left(\hat{\rho}^n\right)
\end{equation}
It can be seen from the above equation that the density perturbations are damped in the HLLE scheme and the result is identical to the linearized analysis.
The x-momentum perturbation equation is
\begin{equation}
\hat{\rho{}u}^{n+1}=\hat{\rho{}u}^n-2\nu\left( \hat{\rho{}u}^n\right)
\end{equation}
The energy perturbation equation is 
\begin{equation}
\dfrac{1}{\gamma-1}\hat{p}^{n+1}+\dfrac{1}{2}\hat{(\rho{}q^2)}^{n+1}=\dfrac{1}{\gamma-1}\hat{p}^n+\dfrac{1}{2}\hat{\rho{}q^2}^n-2\nu \left(\dfrac{1}{\gamma-1}\hat{p}^n+\dfrac{1}{2}\hat{\rho{}q^2}^n\right)
\end{equation}
Expanding the terms and cancelling the terms of the continuity and momentum equations, we obtain the following relation for the pressure perturbation
\begin{equation}
\hat{p}^{n+1}=\hat{p}^n-2\nu\hat{p}^n
\end{equation}
The deviations of the state variables are obtained as
\begin{equation} \label{hll-d1}
\hat{\rho}^{n+1}=\hat{\rho}^n-2\nu\hat{\rho}^n \hspace{1cm} \hat{\rho{}u}^{n+1}=\hat{\rho{}u}^n-2\nu\hat{\rho{}u}^n  \hspace{1cm} \hat{p}^{n+1}=\hat{p}^n-2\nu\hat{p}^n  
\end{equation} 
The terms required in equation (\ref{lf1}) for evaluating the stability of the scheme can be written as 
\begin{equation}
\begin{array}{l}
\hat{\rho}^{2,n+1}-\hat{\rho}^{2,n}=-4\nu(1-\nu)\hat{\rho}^{2,n} \\
\hat{\rho{}u}^{2,n+1}-\hat{\rho{}u}^{2,n}=-4\nu(1-\nu)\hat{\rho{}u}^{2,n}  \\
\hat{p}^{2,n+1}-\hat{p}^{2,n}=-4\nu(1-\nu)\hat{p}^{2,n}
\end{array}
\end{equation}
\begin{equation}
V(\mathbf{U})^{n+1}-V(\mathbf{U})^n=-\dfrac{a_0^2}{\rho{}_0}4\nu(1-\nu)\hat{\rho}^{2,n}-\dfrac{a_0^2}{\rho_0u_0^2}4\nu(1-\nu)\hat{\rho{}u}^{2,n}-\dfrac{1}{\rho_0a_0^2}4\nu(1-\nu)\hat{p}^{2,n}
\end{equation}
Here, the CFL number $\nu$ is always less than unity and $\hat{\rho}^{2,n} >0$, $\hat{\rho{}u}^{2,n}>0$, $\hat{p}^{2,n}>0$. Therefore, 
\begin{equation}
V(\mathbf{U})^{n+1}-V(\mathbf{U})^n< 0 \hspace{1cm} \text{for all } \mathbf{U}^n\neq 0 \hspace{0.3cm} \text{in} \hspace{0.3cm} \mathcal{D}
\end{equation}
It is also observed that  $V(\mathbf{U})\to \infty$ as $||\mathbf{U}|| \to \infty$ for the chosen quadratic Lyapunov function. Therefore, the HLLE Scheme is globally, asymptotically stable. 
\subsection{Analysis of the Roe, HLLEM and HLLC schemes}
For the HLLEM and Roe schemes, the deviation of the state variables are obtained  as
\begin{equation} \label{roe-d1}
\hat{\rho}^{n+1}=\hat{\rho}^n-2\nu\dfrac{\hat{p}^n}{a_0^2} \hspace{1cm} \hat{\rho{}u}^{n+1}=\hat{\rho{}u}^n-2\nu{}u_0\dfrac{\hat{p}^n}{a_0^2}  \hspace{1cm} \hat{p}^{n+1}=\hat{p}^n-2\nu\hat{p}^n  
\end{equation} 
Therefore, the terms required in equation (\ref{lf1}) for evaluating the stability of the scheme are 
\begin{equation}
\begin{array}{l}
\hat{\rho}^{2,n+1}-\hat{\rho}^{2,n}=-4\nu\hat{\rho}^n\dfrac{\hat{p}^n}{a_0^2} +4\nu^2\dfrac{\hat{p}^{2,n}}{a_0^4} \\ 
\hat{\rho{}u}^{2,n+1}-\hat{\rho{}u}^{2,n}=\dfrac{-4\nu{}u_0}{a_0^2}\hat{\rho{}u}^n\hat{p}^n+\dfrac{4\nu^2u_0^2}{a_0^4}\hat{p}^{n,2} \\
\hat{p}^{2,n+1}-\hat{p}^{2,n}=-4\nu(1-\nu)\hat{p}^{2,n}
\end{array}
\end{equation}
\begin{equation}\label{lf-hllem-roe}
\Delta{}V(\mathbf{U}):=V(\mathbf{U})^{n+1}-V(U)^n=-\dfrac{4\nu}{\rho_0}\hat{\rho}^n\hat{p}^n+\dfrac{4\nu^2}{\rho_0a_0^2}\hat{p}^{2,n}-\dfrac{4\nu}{\rho_0u_0}\hat{\rho{}u}^n\hat{p}^n+\dfrac{4\nu^2}{\rho_0a_0^2}\hat{p}^{2,n}-\dfrac{1}{\rho_0a_0^2}4\nu(1-\nu)\hat{p}^{2,n}
\end{equation}
From the above equation, we observe that the delta change in the Lyapunov function can be greater than zero if the deviations in $\rho$ and $\rho{}u$ have opposite signs with respect to the deviation in pressure $p$. Therefore, $
V(\mathbf{U})^{n+1}-V(\mathbf{U})^n\nleq 0 \hspace{2mm} \forall \hspace{3mm} \mathbf{U}^n\neq 0 \hspace{0.3cm} \text{in} \hspace{0.3cm} \mathcal{D}$ and hence, the Roe, HLLEM and HLLC schemes can be unstable. It may be noted that Quirk \cite{quirk} had observed that the Roe scheme may become unstable if the pressure and density perturbations have opposite signs. In the present work, it is observed that the Roe scheme may become unstable if the perturbations in density $\rho$ and transverse momentum $\rho{}u$ have opposite signs with respect to the pressure perturbations. Thus, additional information about the instability mechanism in the Roe, HLLEM and HLLC schemes is obtained here from the direct Lyapunov method.

\begin{figure}[H]
	\begin{center}
	\includegraphics[width=100pt]{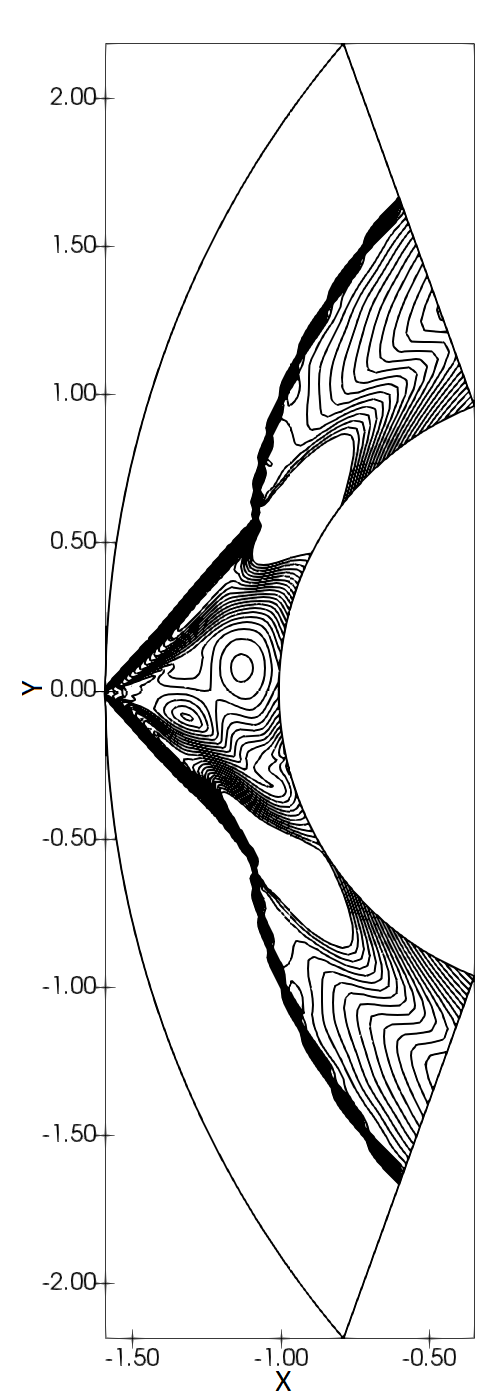} \hspace{1cm} \includegraphics[width=100pt]{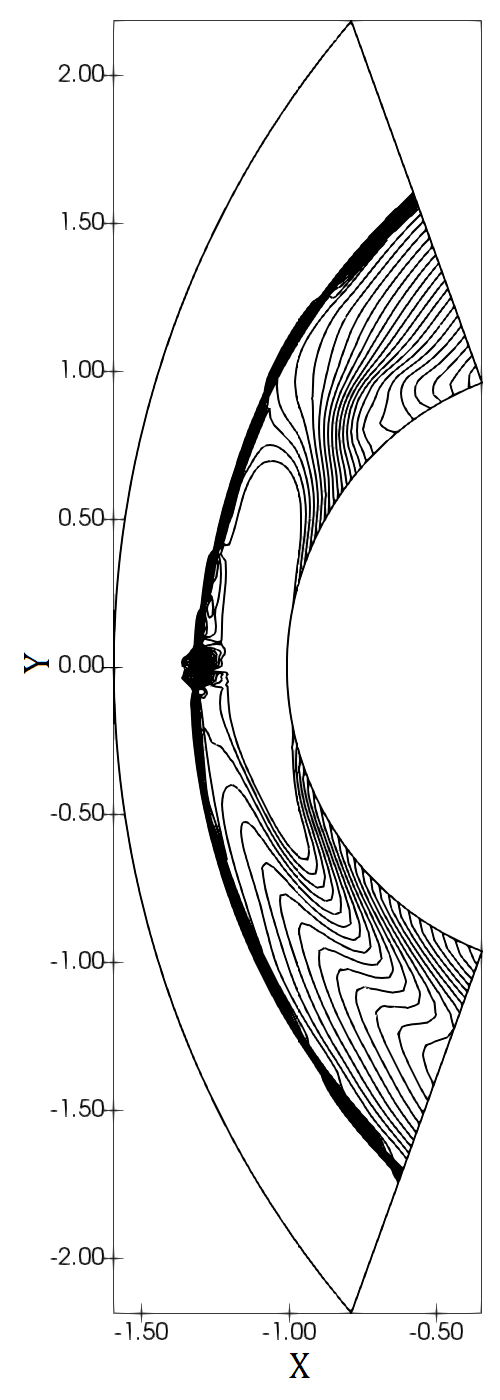} \\
		(a) HLLEM Scheme \hspace{15 mm} (b) HLLEM-LM Scheme 
		\caption{Density contours for $M_{\infty}=20$ flow over a blunt body computed by the HLLEM and HLLEM-LM schemes}
		\label{hllem-lm-comp}
	\end{center}
\end{figure}

Further, it is observed that, for the Roe, HLLEM and HLLC schemes
\begin{equation}
\Delta{}V(\mathbf{U}):=V(\mathbf{U})^{n+1}-V(\mathbf{U})^n= 0 \hspace {1cm} \text{if} \hspace{1cm} \hat{p}=0
\end{equation}
This implies that the Roe, HLLEM and HLLC schemes are simply stable in the absence of pressure perturbation. In other words, it is the pressure perturbation that drives the Roe, HLLEM and HLLC schemes toward numerical instability. Further, equation (\ref{lf-hllem-roe}) also implies that the magnitude of the numerical shock instabilities depends upon the magnitude of the pressure perturbations and hence the numerical instabilities in the HLLEM/Roe schemes can be minimized if the pressure perturbations are minimized. Here, it is worth noting that, Guillard and Viozat \cite{guillard1} and Rieper \cite{rieper2} have shown that in the limit as Mach number tends to zero, the pressure fluctuations in the Roe scheme for the discrete Euler equations are of the order of Mach number $\mathcal{O} (M_*)$, while for the continuous Euler equations, the pressure fluctuations are of the order of Mach number squared $\mathcal{O} (M_*^2)$. Guillard and Murrone \cite{guillard2} have studied the Godunov scheme and have shown that the interface pressure in the Riemann problem contains acoustic waves of order $\mathcal{O} (M_*)$ even if the initial data are well prepared and contain pressure fluctuations of order $\mathcal{O}(M_*^2)$. Several low Mach corrections have been proposed to the approximate Riemann solvers like the Roe scheme such that the pressure fluctuations in these schemes are consistent with the continuous Euler equations. With the low Mach corrections of Li and Gu \cite{li-roe}, Rieper \cite{rieper2}, Thornber et al. \cite{th1}, and Dellacherie et al. \cite{della2}, the pressure perturbations in the approximate Riemann solvers at low Mach numbers can be pushed to the order of Mach number squared. Therefore, we conjecture that with appropriate low Mach corrections, the pressure perturbations in the plane transverse to a grid-aligned normal shock can be reduced to the order of Mach number squared and the numerical shock instabilities can also be correspondingly reduced. To verify the conjecture, computations are carried out for the Mach 20 flow around a blunt body, which is the standard test case for the carbuncle phenomenon, using the HLLEM scheme and the HLLEM scheme with Dellacherie et al.'s low Mach fix. The computational details can be found in \ref{carbuncle}. The comparison of the results of the original HLLEM scheme and the HLLEM scheme with Dellacherie et al.'s fix is shown in Fig. \ref{hllem-lm-comp} for a Mach 20 flow around a blunt body. It can be seen from the figure that the magnitude of the shock instability has reduced significantly in the low Mach version of the HLLEM scheme. However, shock instability is still present in the low Mach version of the HLLEM scheme as the pressure perturbations, which are reduced in magnitude but still being generated, feed the perturbations in density and transverse momentum. The results validate the observation that  the magnitude of numerical shock instabilities in the complete or all-wave Riemann solvers depend upon the magnitude of pressure perturbations and the magnitude of the numerical shock instabilities can be reduced by reducing the magnitude of pressure perturbations with the help of low Mach corrections.
\subsection{Analysis of the HLLCM and HLLEC schemes}
For the HLLCM and HLLEC schemes, the deviation of the state variables are
\begin{equation} \label{hllcps-d1}
\hat{\rho}^{n+1}=\hat{\rho}^n-\dfrac{2\nu\hat{p}^n}{a_0^2} \hspace{1cm} \hat{\rho{}u}^{n+1}=\hat{\rho{}u}^n-2\nu{}\left(\hat{\rho{}u}^n-u_0\hat{\rho}+u_0\dfrac{\hat{p}^n}{a_0^2} \right) \hspace{1cm} \hat{p}^{n+1}=\hat{p}^n-2\nu\hat{p}^n  
\end{equation} 
The expression for the delta change in the Lyapunov function for the HLLCM/HLLEC schemes becomes complicated. Hence, we evaluate the delta change in the Lyapunov function for a special case of zero shear velocity perturbation, i.e., $\hat{u}=0$. For the special case of zero shear velocity perturbation, $\hat{\rho{}u}=u_o\hat{\rho}$. Therefore, the terms required in equation (\ref{lf1}) for evaluating the stability of the scheme are
\begin{equation}
\begin{array}{l}
\hat{\rho}^{2,n+1}-\hat{\rho}^{2,n}=-4\nu\hat{\rho}^n\dfrac{\hat{p}^n}{a_0^2} +4\nu^2\dfrac{\hat{p}^{2,n}}{a_0^4} \\
\hat{\rho{}u}^{2,n+1}-\hat{\rho{}u}^{2,n}=-4\nu{} \dfrac{u_0}{a_0^2}\hat{\rho{}u}^n\hat{p}^n+ 4\nu^2\dfrac{u_0^2}{a_0^4}\hat{p}^{2,n}\\
\hat{p}^{2,n+1}-\hat{p}^{2,n}=-4\nu(1-\nu)\hat{p}^{2,n}
\end{array}
\end{equation}
\begin{equation}
\Delta{}V(\mathbf{U}):=V(\mathbf{U})^{n+1}-V(U)^n=-\dfrac{4\nu}{\rho_0}\hat{\rho}^n\hat{p}^n+\dfrac{4\nu^2}{\rho_0a_0^2}\hat{p}^{2,n}-\dfrac{4\nu}{\rho_0u_0}\hat{\rho{}u}^n\hat{p}^n+\dfrac{4\nu^2}{\rho_0a_0^2}\hat{p}^{2,n}-\dfrac{1}{\rho_0a_0^2}4\nu(1-\nu)\hat{p}^{2,n}
\end{equation}
From the above equation, we observe that the delta change in the Lyapunov function can be greater than zero under the following conditions
\begin{itemize}
\item  $\hat{\rho}^n$ and $\hat{p}^n$ have opposite signs
\item $\hat{\rho{}u}^n$ and $\hat{p}^n$ have opposite signs
\end{itemize}
Therefore, $ V(\mathbf{U})^{n+1}-V(\mathbf{U})^n\nleq 0 \hspace{1mm} \forall \hspace{1mm} \mathbf{U}^n\neq 0 \hspace{0.3cm} \text{in} \hspace{0.3cm} \mathcal{D}$ for the HLLCM and HLLEC scheme. Hence, the HLLCM and HLLEC scheme are not always stable. The mechanism of numerical shock instability in the HLLCM and HLLEC schemes is similar to the Roe/HLLEM scheme and the pressure perturbations feeding the density and transverse momentum perturbations are the cause of instability in the HLLCM/HLLEC scheme. It may be noted that the shear velocity perturbations are damped in the HLLCM/HLLEC schemes, unlike the Roe/HLLEM schemes. It is conjectured that the stability of the HLLCM/HLLEC schemes may improve if non-zero shear velocity perturbations are considered. 
\subsection{Analysis of the HLLS and HLLES schemes}
For the HLLS and HLLES schemes, the deviation of the state variables are
\begin{equation} \label{hlles-d1}
\hat{\rho}^{n+1}=\hat{\rho}^n-2\nu\hat{\rho}^n \hspace{1cm} \hat{\rho{}u}^{n+1}=\hat{\rho{}u}^n-2\nu{}u_0\hat{\rho{}}^n  \hspace{1cm} \hat{p}^{n+1}=\hat{p}^n-2\nu\hat{p}^n  
\end{equation} 
The terms required in equation (\ref{lf1}) for evaluating the stability of the scheme are
\begin{equation}
\begin{array}{l}
\hat{\rho}^{2,n+1}-\hat{\rho}^{2,n}=-4\nu(1-\nu)\hat{\rho}^{2,n} \\
\hat{\rho{}u}^{2,n+1}-\hat{\rho{}u}^{2,n}=-4\nu{}u_0\hat{\rho}\hat{\rho{}u}^n +4\nu^2u_0^2\hat{\rho}^2 \\
\hat{p}^{2,n+1}-\hat{p}^{2,n}=-4\nu(1-\nu)\hat{p}^{2,n}
\end{array}
\end{equation}
\begin{equation}
\Delta{}V(\mathbf{U}):=V(\mathbf{U})^{n+1}-V(\mathbf{U})^n=-\dfrac{a_0^2}{\rho{}_0}4\nu(1-2\nu)\hat{\rho}^{2,n}-\dfrac{a_0^2}{\rho_0u_0}4\nu\hat{\rho}\hat{\rho{}u}^n-\dfrac{1}{\rho_0a_0^2}4\nu(1-\nu)\hat{p}^{2,n}
\end{equation}
Here, the CFL number $\nu$ is always less than unity and the terms $\hat{\rho}^{2,n}$, and $\hat{p}^{2,n}$ are always positive. It can be seen from the equation that the delta change in the Lyapunov function can become positive if $\hat{\rho}$ and $\hat{\rho{}u}/u_0$ have opposite signs and $\hat{\rho{}u}/u_0>>\hat{\rho}$. Therefore, for the HLLES scheme $V(\mathbf{U})^{n+1}-V(\mathbf{U})^n\nleq 0 \hspace{1mm} \forall \hspace{1mm} \mathbf{U}^n\neq 0 \hspace{0.3cm} \text{in} \hspace{0.3cm} \mathcal{D} $. Hence, the HLLS and HLLES schemes are not always stable. Further, it is observed that, for the HLLS/HLLES schemes
\begin{equation}
\Delta{}V(\mathbf{U}):=V(\mathbf{U})^{n+1}-V(\mathbf{U})^n \le 0 \hspace {1cm} \text{if} \hspace{1cm} \hat{\rho}=0
\end{equation}
This implies that the HLLS/HLLES schemes are stable in the absence of density perturbation and density perturbation feeding the perturbation in transverse momentum $\rho{}u$ is responsible for instability in the HLLS/HLLES scheme.

\section{A formulation for a shock-stable complete approximate Riemann solver}
Stability analysis of the complete approximate Riemann solvers like Roe, HLLEM and HLLC schemes using the direct Lyapunov method indicates that the pressure perturbations feeding the perturbations in density and transverse momentum $\rho{}u$ are responsible for the numerical shock instabilities in the scheme. Therefore, the instabilities in these complete Riemann solvers can be suppressed by reducing the rate at which the pressure perturbations feed the density and transverse momentum perturbations and by introducing damping of the density and transverse momentum perturbations in the presence of pressure perturbations. The stability analysis and numerical results also indicate that if the magnitude of the pressure perturbations is reduced, the magnitude of the numerical shock instabilities will also reduce. It has been shown using asymptotic analysis that in the limit of zero Mach number the magnitude of the pressure perturbations in the Roe scheme are of the order of Mach number \cite{guillard1, rieper2}. In case low Mach corrections like those proposed by Rieper \cite{rieper2}, Thornber et al. \cite{th1} or Dellacherie et al. \cite{della2} are applied to these schemes, the pressure perturbations in the plane transverse to the shock shall be of the order of Mach number squared, like the continuous Euler equations. Hence, the magnitude of the pressure perturbations shall be lower in case low Mach corrections are applied to these scheme and the magnitude of the numerical shock instabilities shall be proportionally lower.

Based on these observations, a shock-stable HLLEM flux function including the low Mach fix of Dellacherie et al. \cite{della2} for Godunov-type schemes, can be written as
\begin{equation} \label{hllem-fp1d-flux}
\mathbf{F}(\mathbf{U}_R,\mathbf{U}_L)=\dfrac{S_R\mathbf{F}_L-S_L\mathbf{F}_R}{S_R-S_L}+\dfrac{S_RS_L}{S_R-S_L}\left(\Delta{}\mathbf{U}-\sum_{k=2}^3\delta_{k,new}\tilde{\alpha_k}\tilde{\mathbf{R}}_k\right)-(1-\theta)\tilde{\rho}\tilde{a}\left[0,\;\Delta{}u_n,\;0,\;0\right]^T
\end{equation} 

where 

$\tilde{\alpha}_2=\Delta{}\rho-\dfrac{\Delta{}p}{\tilde{a}^2}$ and $\tilde{\alpha}_3=\tilde{\rho}\Delta{}u_t$ are the wave strengths of the contact and shear waves respectively,

$\tilde{\mathbf{R}}_2=\left[1,\;\tilde{u}_n,\;\tilde{u}_t,\;\frac{1}{2}(\tilde{u}_n^2+\tilde{u}_t^2)\right]^T$ and $\tilde{\mathbf{R}}_3=\left[0,\; 0,\; 1,\;\tilde{u}_t\right]^T$ are the right eigenvectors of the contact and shear waves respectively in the flux Jacobian matrix, 

$\delta_{2,new}$ and $\delta_{3,new}$ are the proposed anti-diffusion coefficients for the contact and shear waves respectively,

$\theta=min\left[max(M_R,M_L),1\right]$ is the local Mach number function of Dellacherie et al. and $\Delta(.)=(.)_R-(.)_L$.
 
In order to reduce the rate of feeding of the pressure perturbations and to increase the rate of damping of density and transverse momentum perturbations in presence of the pressure perturbations, modifications in the anti-diffusion coefficients are proposed. The proposed anti-diffusion coefficients are defined as 
\begin{equation}\label{delta_1d}
\delta_{k,new}=\delta_k\left(1-\left(\dfrac{\Delta{}p}{p_{max}}\right)^r \right)
\end{equation}

where $\delta_2=\delta_3=\frac{\tilde{a}}{\tilde{a}+|\tilde{u_n}|}$ are the original anti-diffusion coefficients for the contact and shear waves respectively \cite{park}, 

$\Delta{}p=|p_L-p_R|$, 

$p_{max}=max(p_L,p_R)$ and 

$r$ is a term that determines the reduction in the feeding rate of the pressure perturbations as well as damping of the density and transverse momentum perturbations. 

Numerical experiments suggest that $r=\frac{1}{3}$ is sufficient for suppressing the numerical shock instabilities in the HLLEM scheme and hence $r=\frac{1}{3}$ is used in the present work. It may be noted here that the anti-diffusion coefficients are modified with a one-dimensional pressure function in the present work. On the other hand, the modification proposed by Simon and Mandal \cite{san-hllem} and Xie et al. \cite{xie-hllem} and several others comprise scaling the anti-diffusion coefficients with a multi-dimensional shock-sensing pressure function which require a much larger stencil. Therefore, the proposed HLLEM flux shown in equation (\ref{hllem-fp1d-flux}) is named as HLLEM-FP1D with FP1D indicating a one-dimensional pressure function.

\begin{figure}[H]
	\begin{center}
		\includegraphics[width=300pt]{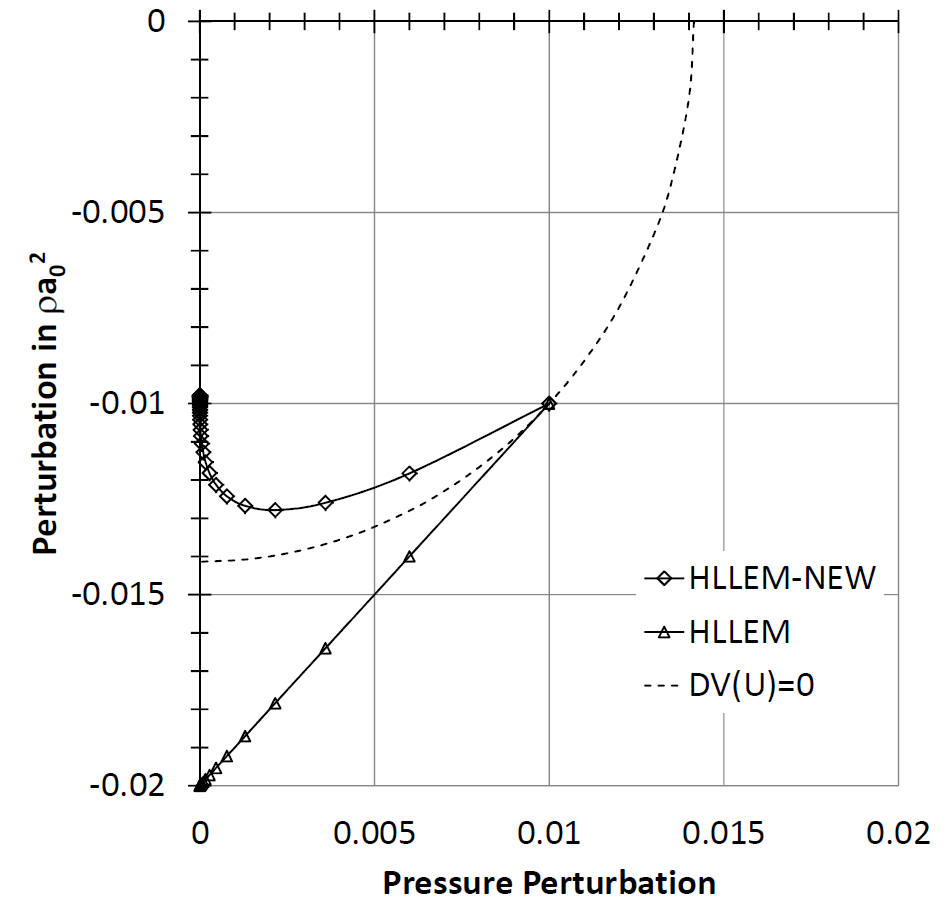} 
		\caption{Phase portrait of the HLLEM and proposed HLLEM-FP1D schemes for positive pressure and negative density  perturbations}
	\label{phase-portrait}
	\end{center}
\end{figure}

Stability analysis of the proposed HLLEM-FP1D scheme is carried out using the direct Lyapunov method shown in Section \ref{direct-lyapunov} to determine its stability characteristics. Since the normal velocity is zero for the case of  grid aligned shock, the anti-diffusion coefficients of the contact and shear waves shown in equation \ref{delta_1d} become
\begin{equation}\label{delta_hllemfp1d}
\delta=\delta_2=\delta_3=\left( 1-\left(\dfrac{2|\hat{p}|^n}{p_0+|\hat{p}|^n}\right)^{1/3}\right)
\end{equation}
The deviations in the state variables in the proposed HLLEM-FP1D scheme are obtained as 
\begin{equation} \label{perturbation_hllemfp1d}
\begin{array}{l}
\hat{\rho}^{n+1}=\hat{\rho}^n-2\nu\left((1-\delta)\hat{\rho}^n+\delta{}\dfrac{\hat{p}^n}{a_0^2} \right) \\
\hat{\rho{}u}^{n+1}=\hat{\rho{}u}^n-2\nu{}\left((1-\delta)\hat{\rho{}u}^n+\delta{}u_0\dfrac{\hat{p}^n}{a_0^2}\right)  \\
\hat{p}^{n+1}=\hat{p}^n-2\nu\hat{p}^n  
\end{array}
\end{equation} 
It can be seen from equations (\ref{roe-d1}, \ref{delta_hllemfp1d}, \ref{perturbation_hllemfp1d}) that in case of a pressure perturbation, damping of the density and transverse momentum  perturbations are introduced  unlike the original HLLEM scheme and the feeding rate of the pressure perturbation is reduced as compared to the original HLLEM scheme. It can also be shown from the equation that the damping rate of the density and transverse momentum perturbations exceed the feeding rate of the pressure perturbations. For example, for a pressure perturbation of $\hat{p}=0.001$, the anti-diffusion coefficient $\delta$ becomes 0.90 and the damping rate $(1-\delta)$ becomes 0.10 and exceeds the feeding rate of $\hat{p}=0.001$.  
It is observed that the anti-diffusion coefficient in the proposed scheme is highly non-linear. Is it difficult to demonstrate the stability of such non-linear system for all conditions. However, the stability of the proposed scheme can however be illustrated with phase portraits. The phase portrait of the HLLEM and the proposed HLLEM-FP1D schemes are shown in Fig. \ref{phase-portrait} for the critical case of an initial positive pressure perturbation and negative density perturbation. It can be seen from the figure that the delta change in the Lyapunov function for the HLLEM scheme is positive $(\Delta{}V(\mathbf{U})>0)$ indicating that the perturbed state has moved away from the steady state point. On the other hand, the delta change in the Lyapunov function for the HLLEM-FP1D scheme is negative $(\Delta{}V(\mathbf{U})<0)$ indicating that the perturbed state has moved closer to the steady state point. The pressure perturbation has been fully damped to zero while the density perturbation in the proposed HLLEM-FP1D scheme has settled to a value less than the initial perturbation. The pressure perturbation is unable to feed the density perturbations in the proposed scheme, unlike the HLLEM scheme. The phase portrait of the HLLEM-FP1D indicates that the scheme is stable and hence the scheme is expected to be free from the numerical shock instability problems of the original HLLEM scheme.
In the next section, numerical results are presented for several test cases to demonstrate the robustness of the proposed scheme to the numerical shock instabilities.

\section{Results}
The performance of the proposed HLLEM-FP1D scheme is evaluated by solving the problems involving planar shock, double Mach reflection, forward-facing step, blunt body, and expansion corner. The performance of the proposed scheme is compared with the HLLE and HLLEM schemes. First-order results are shown for the test cases since the presence of numerical shock instabilities is prominently visible in first-order schemes. Second-order results are also presented for a few test cases to demonstrate the ability of the proposed scheme to resolve the flow features.
\subsection{Planar shock problem}
The problem consists of a Mach 6.0 shock wave propagating down in a rectangular channel. The domain consists of 800$\times{}$20 cells. The centreline in the y direction (11$^{th}$ grid line) is perturbed to promote odd-even decoupling along the length of the shock as $y_{i,11}=\left\lbrace\begin{array}{c}\ y_{i,11}+0.001 \hspace{4mm}\text{if  i is even}\\ 	y_{i,11}-0.001 \hspace{4mm}\text{if i is odd}\end{array}\right.$. 

\begin{figure}[H]
	\begin{center}
	\includegraphics[width=225pt]{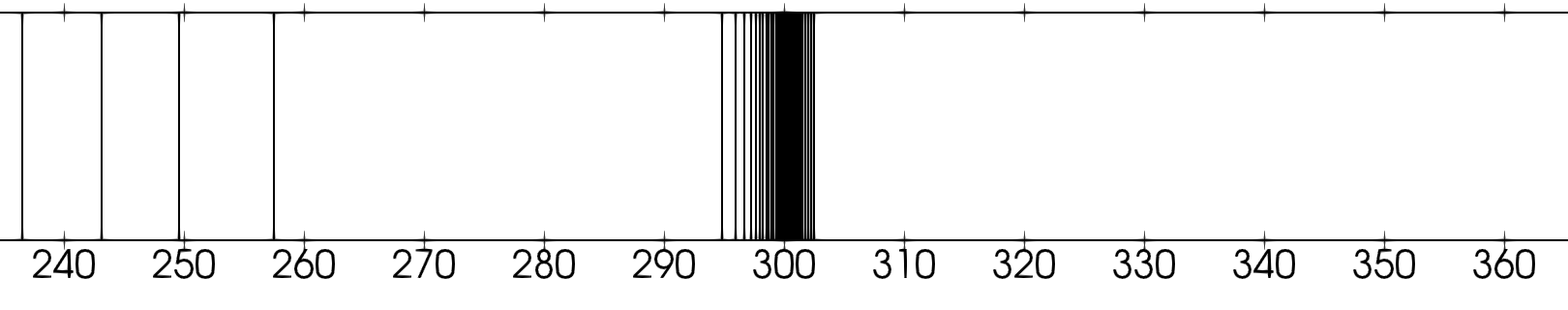} \includegraphics[width=225pt]{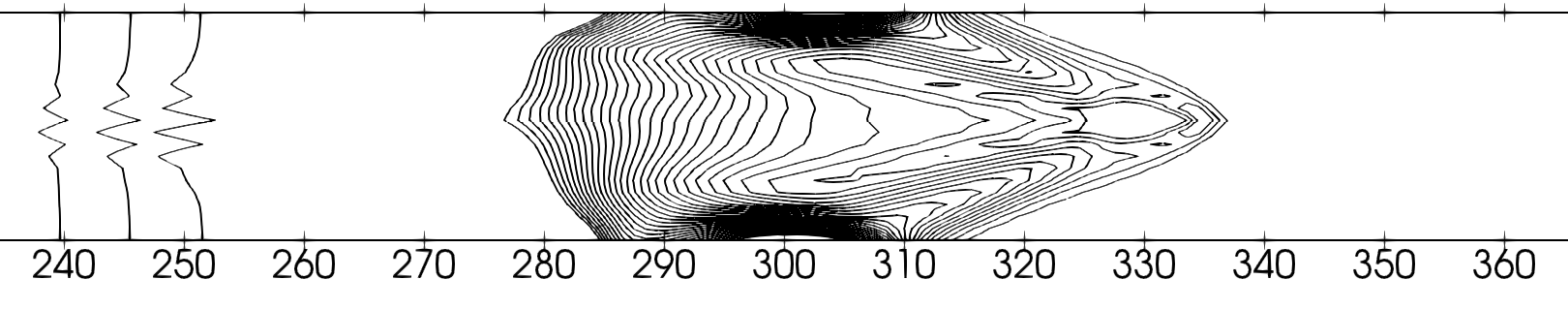} \\ 
	(a) HLLE Scheme \hspace{5cm} (b) HLLEM Scheme\\
	\includegraphics[width=225pt]{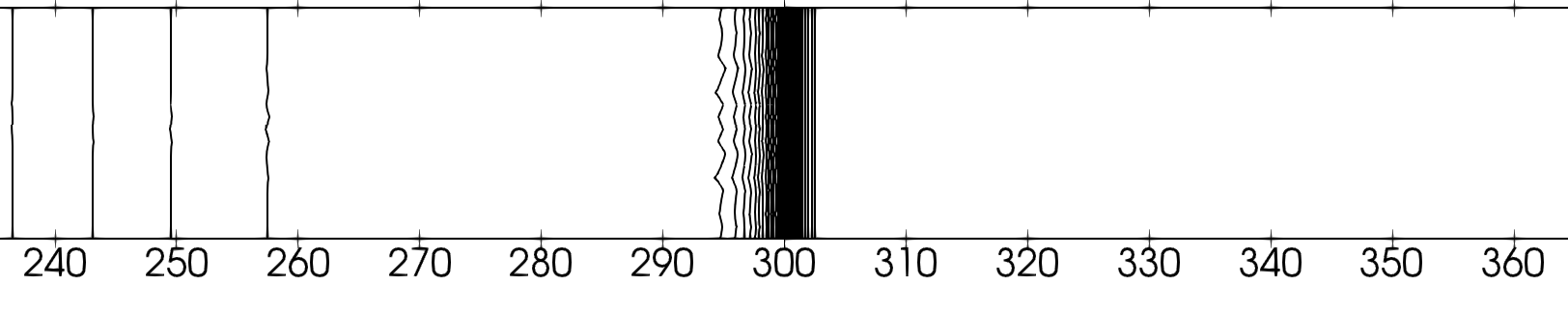} \includegraphics[width=225pt]{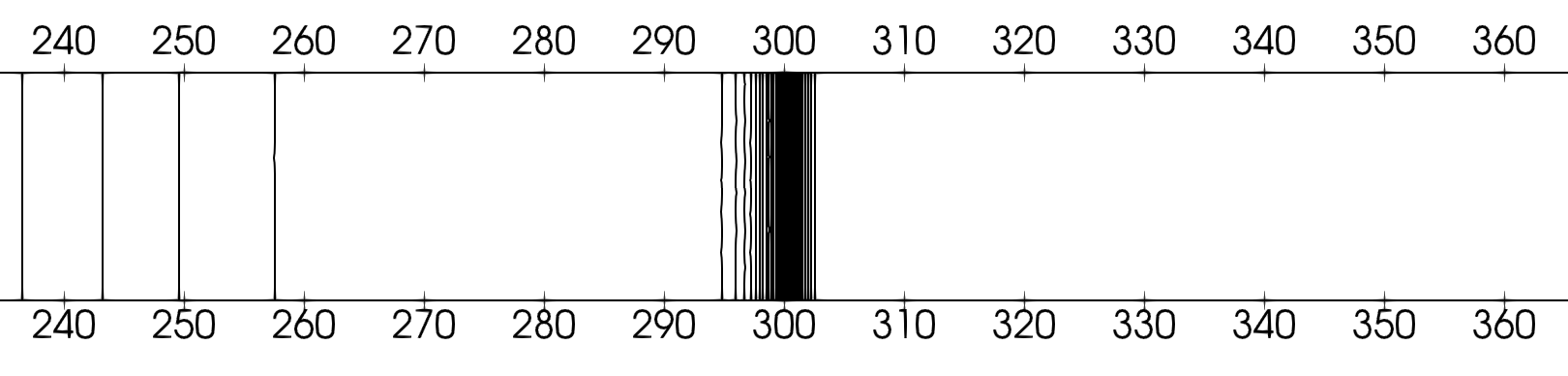}\\ 
	(c)  HLLEM-FP1D Scheme $r=\frac{1}{2}$ \hspace{2cm} (d)  HLLEM-FP1D Scheme $r=\frac{1}{3}$\\ 
	\caption{Density contours for the $M_{\infty}=6$ planar shock problem computed by the HLLE, HLLEM and HLLEM-FP1D schemes. The results are shown at time $t=55$ units.}
	\label{planar-shock-hllem1d}
	\end{center}
\end{figure}

The domain is initialized with  $\rho=1.4$, $p=1.0$, $u=0$, $v=0$. The post-shock values are imposed at the inlet and zero gradient boundary conditions are imposed at the exit. Solid wall boundary conditions are imposed at the top and bottom. The density contour plot for the original and proposed HLLEM schemes at time $t=55$ is shown in Fig. \ref{planar-shock-hllem1d} and 30 contour levels ranging from 1.6 to 7.0 are shown. It can be seen from the figure that the original HLLEM scheme show odd-even decoupling while the proposed HLLEM-FP1D scheme is free from the odd-even decoupling problem.
\subsection{Double Mach reflection problem}
\begin{figure}[H]
	\begin{center}
		\includegraphics[width=225pt]{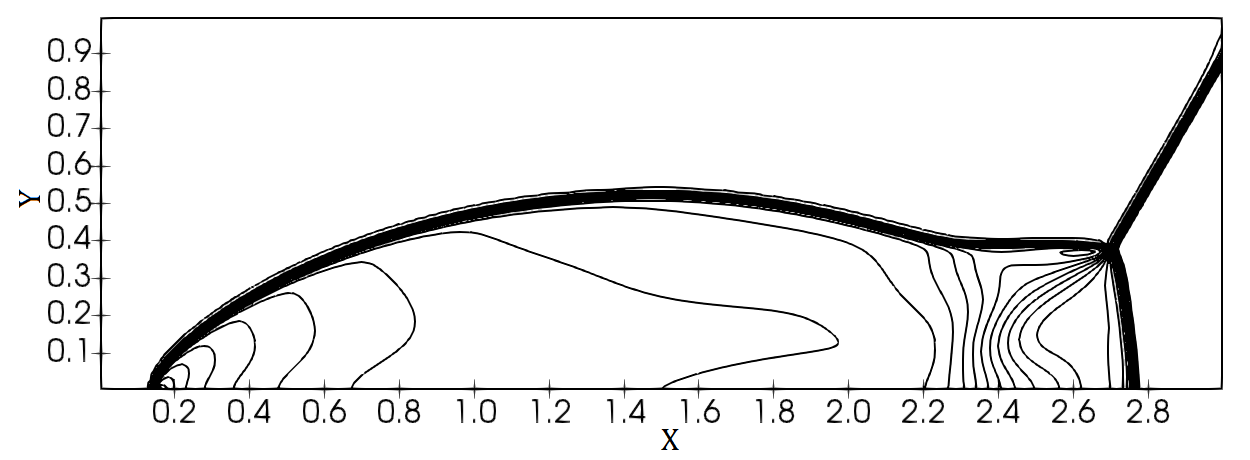} \includegraphics[width=225pt]{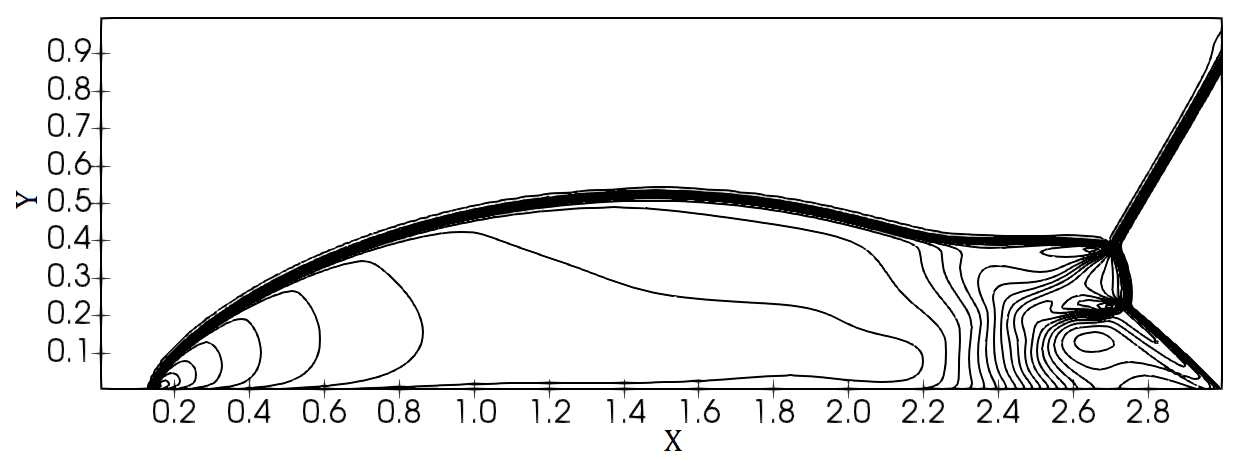}\\ 
		(a) HLLE Scheme \hspace{5cm} (b) HLLEM Scheme \\
		\includegraphics[width=225pt]{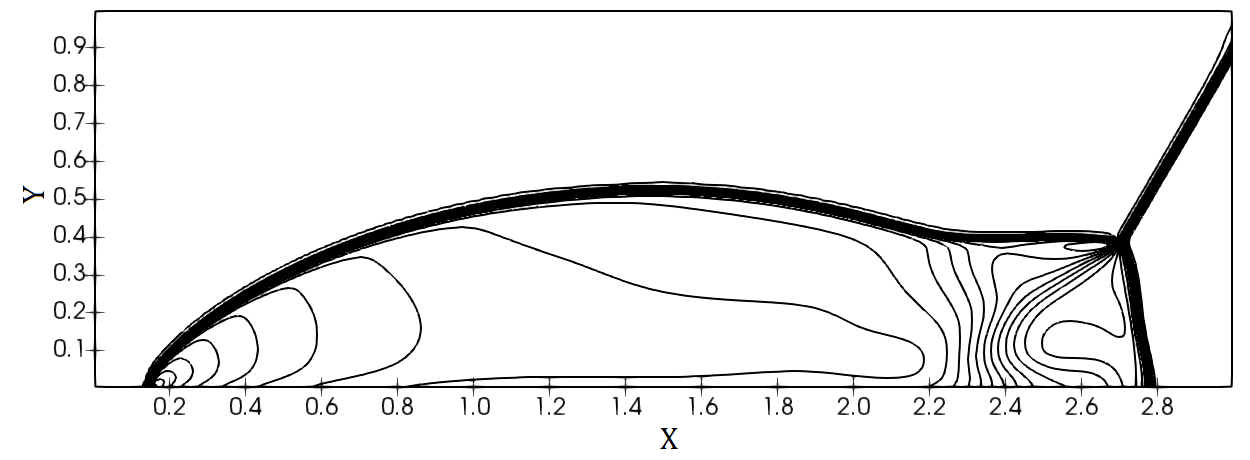} \\ 
		(c) HLLEM-FP1D Scheme $r=\frac{1}{3}$ 
		\caption{Density contours for the $M_{\infty}=10$ double Mach reflection problem computed by the first-order HLLE, HLLEM and HLLEM-FP1D schemes. The results are shown at time $t=0.020026$ units.}
		\label{dmr-hllem1d}
	\end{center}
\end{figure}

\begin{figure}[H]
	\begin{center}
		\includegraphics[width=230pt]{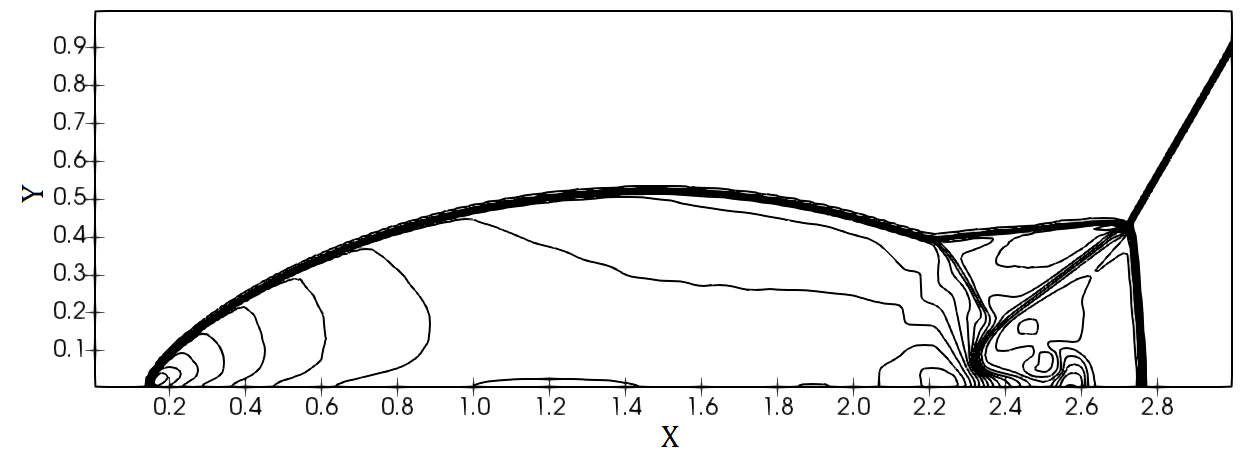} \includegraphics[width=230pt]{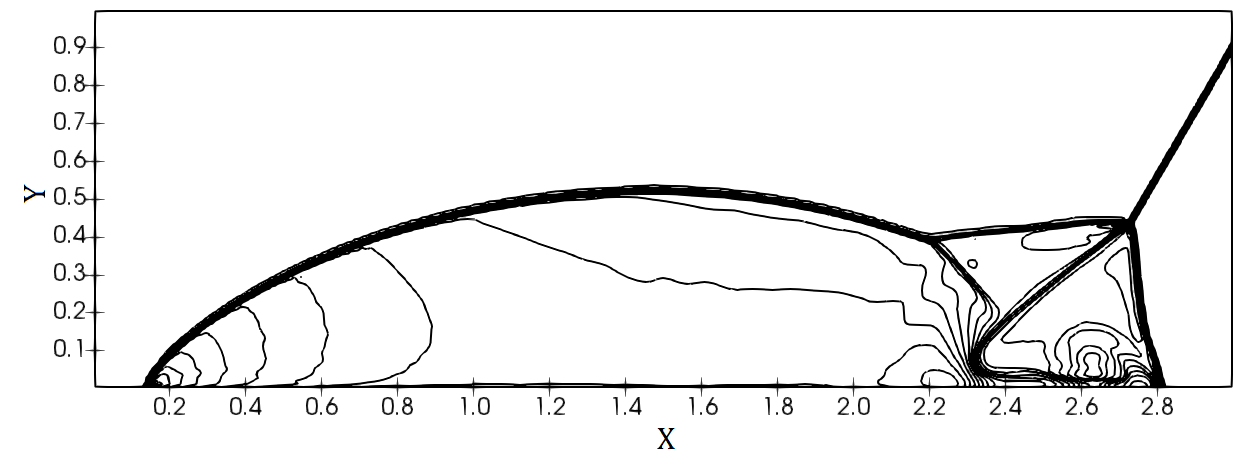} \\
		(a) HLLE Scheme \hspace{5cm} (b) HLLEM Scheme \\
		\includegraphics[width=230pt]{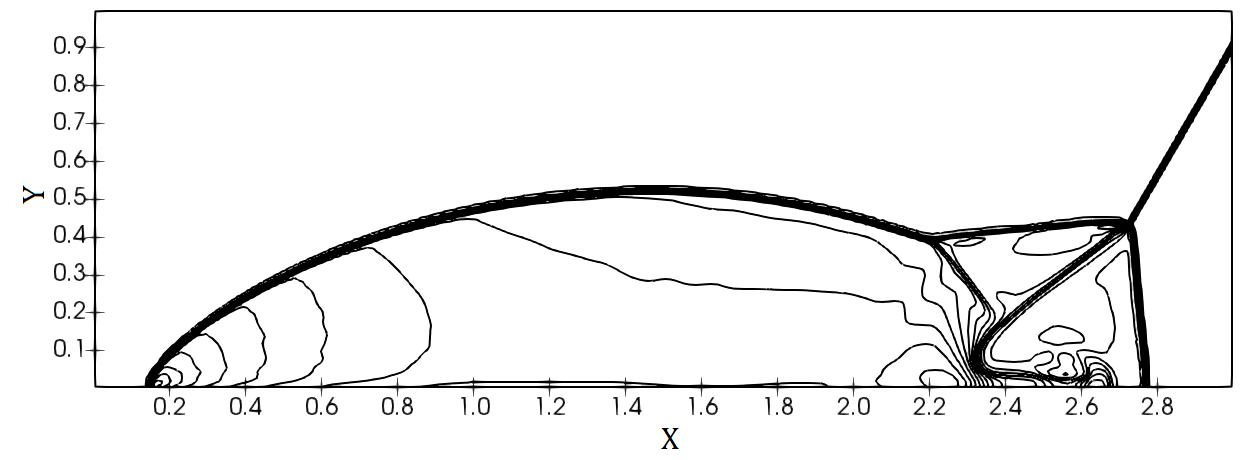}\\ 
		(c) HLLEM-FP1D Scheme
		\caption{Density contours for the $M_{\infty}=10$ double Mach reflection problem computed by the second-order HLLE and HLLEM-FP1D schemes. The results are shown at time $t=0.020026$ units.}
		\label{dmr-hllem1d-o2}
	\end{center}
\end{figure}

The domain is four units long and one unit wide. The domain is divided into 480$\times{}$120 cells. A Mach 10 shock initially making a 60-degree angle at $x=1/6$ with the bottom reflective wall is made to propagate through the domain. The domain ahead of the shock is initialized to pre-shock value ($\rho=1.4$, $p=1$, $u=0$, $v=0$) and domain behind the shock is assigned post-shock values. The inlet boundary condition is set to post-shock values and the outlet boundary is set to zero gradients. The top boundary is set to simulate actual shock movement. At the bottom, the post-shock boundary condition is set up to $x=1/6$ and the reflective wall boundary condition is set thereafter. Density contours are shown in Fig. \ref{dmr-hllem1d} for the original and modified HLLEM schemes at time $t=2.00260\times{}10^{-1}$ units. Thirty (30) contour levels ranging from 2.0 to 21.5 are shown. A severe kinked Mach stem is observed in the original HLLEM scheme, while the proposed HLLEM-FP1D scheme is free from kinked the Mach stem problem. The second-order results of the HLLE and HLLEM-FP1D schemes are shown in Fig. \ref{dmr-hllem1d-o2}. The second-order computations are carried out using MUSCL reconstruction \cite{muscl} with the van Leer limiter. The second-order SSPRK method of Gottlieb and Shu \cite{gott} is used for time integration. The shocks are much thinner and the flow features are better resolved in the second-order results. The second-order results also show that the proposed HLLEM-FP1D scheme is free from the kinked Mach step problem. The second-order results of the HLLE-TNP scheme are comparable to the fifth-order WENO results of Wang et al.\cite{wang} and Zhao et al. \cite{weno6-zhao} who have used a larger grid size of $800\times200$.
\subsection{Forward-facing step problem}
\begin{figure}[H]
	\begin{center}
		\includegraphics[width=225pt]{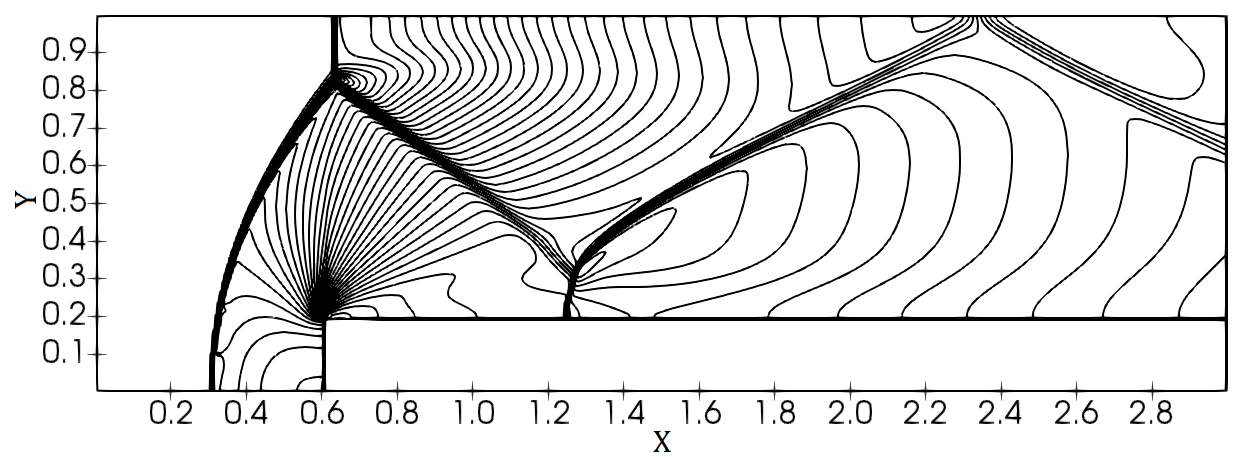} \includegraphics[width=225pt]{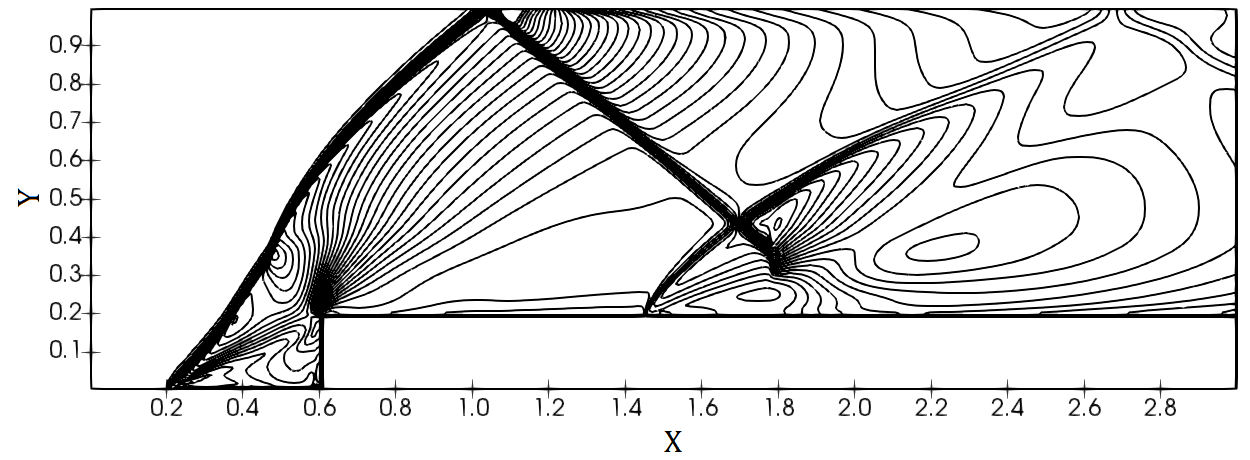}\\
		(a) HLLE Scheme \hspace{5cm} (b) HLLEM Scheme \\
		\includegraphics[width=225pt]{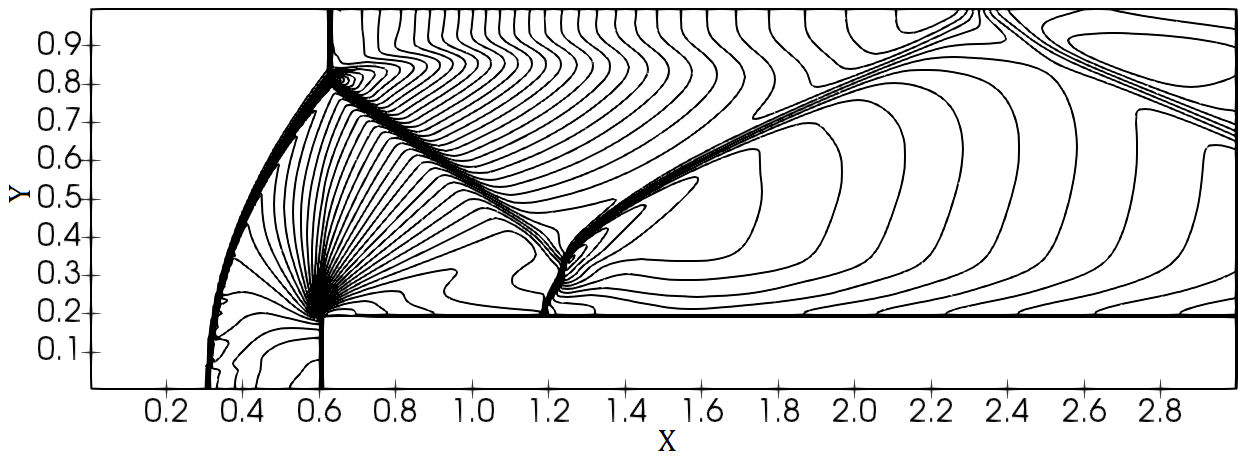} \\
		(c) HLLEM-FP1D Scheme $r=\frac{1}{3}$
		\caption{Density contours for the $M_{\infty}=3$ flow over a forward-facing step computed by the first-order HLLE, HLLEM and HLLEM-FP1D schemes with a grid of $120\times40$. The results are shown at time $t=4$ units.}
		\label{ffs-hllem1d}
	\end{center}
\end{figure}

\begin{figure}[H]
	\begin{center}
		\includegraphics[width=230pt]{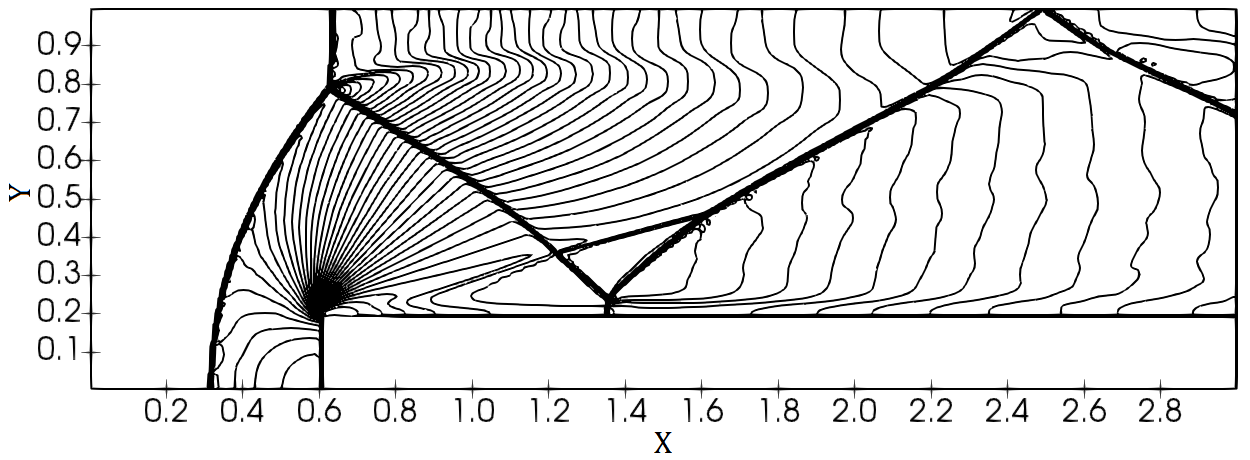} \includegraphics[width=230pt]{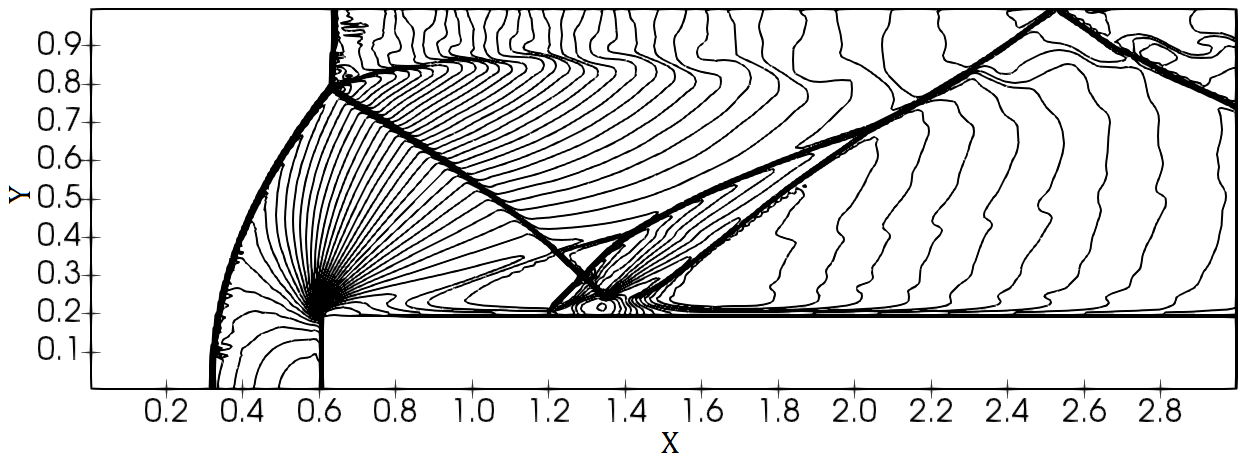} \\
		(a) HLLE Scheme \hspace{5cm} (b) HLLEM Scheme
		\includegraphics[width=230pt]{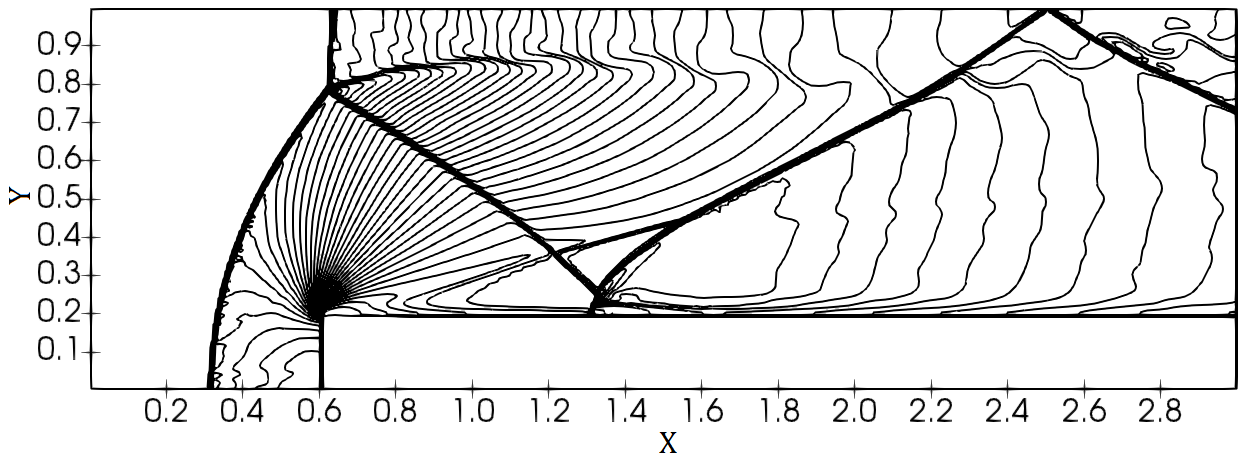} \\
		 (c) HLLEM-FP1D Scheme $r=\frac{1}{3}$
		\caption{Density contours for the $M_{\infty}=3$ flow over a forward-facing step computed by the second-order HLLE and HLLEM-FP1D schemes with a grid of $480\times160$. The results are shown at time $t=4$ units.}
		\label{ffs-hllem1d-o2}
	\end{center}
\end{figure}

The problem consists of a Mach 3 flow over a forward-facing step. The geometry consists of a step located 0.6 units downstream of the inlet and 0.2 units high. A mesh of 480$\times{}$160 cells is used for a domain of 3 units long and 1 unit high. The complete domain is initialized with $\rho=1.4$, $p=1$, $u=3$, and $v=0$. The inlet boundary is set to free-stream conditions, while the outlet boundary has zero gradients. At the top and bottom, reflective wall boundary conditions are set. Density contour plots for the original and modified HLLEM schemes are shown in Fig. \ref{ffs-hllem1d} at time $t=4.0$ units. A total of 45 contours from 0.2 to 7.0 are shown in the figure. A severe carbuncle in the primary bow shock is found to be present in the original HLLEM scheme while the proposed HLLEM scheme captures the primary and reflected shock without numerical instabilities. It is felt that the flow features of the forward-facing step are not well resolved by the first-order computations shown in Fig. \ref{ffs-hllem1d}. Therefore, second-order computations are also performed. The second-order computations are carried out using MUSCL reconstruction \cite{muscl} with the van Leer limiter. The second-order SSPRK method of Gottlieb and Shu \cite{gott} is used for time integration. The second-order results of the HLLE and HLLEM-FP1D schemes are shown in Fig. \ref{ffs-hllem1d-o2}. The contours of the HLLEM-FP1D scheme show slip-lines (contact discontinuity) originating from the triple point near the top wall. The shock are much thinner than the first-order results. Both the classical HLLE and the proposed HLLEM-FP1D schemes resolve the shock structure without encountering numerical instability problems. The slip lines (contact discontinuity) that originate from the Mach stem near the upper wall are found to dissipate rapidly in the HLLE scheme. On the other hand, the slip lines (contact discontinuity) originating from the Mach stem near the upper wall are well captured by the proposed HLLEM-FP1D scheme and the slip lines extend fully downstream of the domain. Further, the phenomenon of roll-up or rolling of the slip lines are also well captured by the HLLEM-FP1D scheme. The density contour plot of the HLLEM-FP1D scheme demonstrate that the scheme is capable of accurate, crisp and stable shock resolution along with contact discontinuity resolution. The second-order results of the HLLEM-FP1D scheme are comparable to the fifth-order WENO results of Zhao et al. \cite{weno6-zhao} who have used a larger grid size of $600\times200$. The HLLEM scheme exhibit mild instability in the bow shock ahead of the step and in the normal shock  near the outer wall. The HLLEM scheme  also exhibit instability in the reflected shocks around the step.  
\subsection{Blunt body problem} \label{carbuncle}
\begin{figure} [H]
	\begin{center}
	\includegraphics[width=100pt]{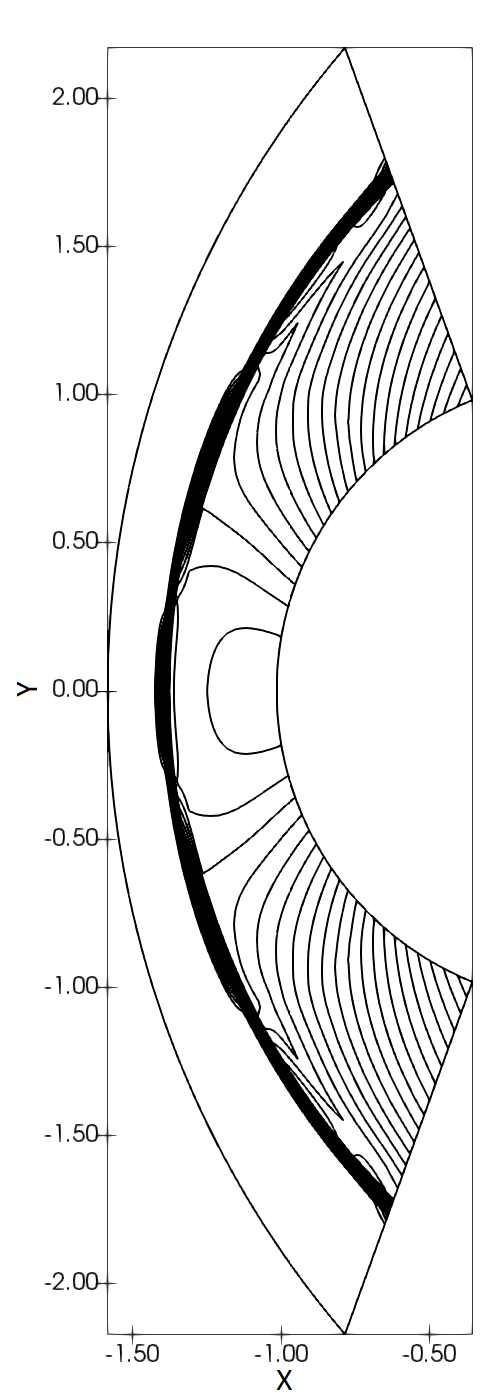} \includegraphics[width=100pt]{carbuncle-hllem-final.png} \includegraphics[width=100pt]{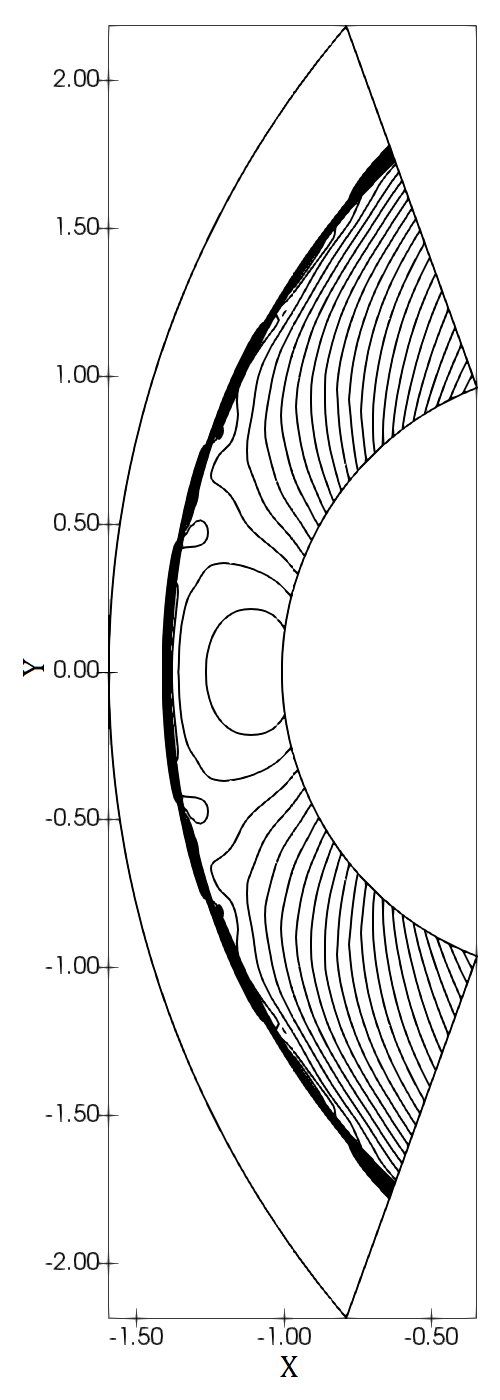}\includegraphics[width=100pt]{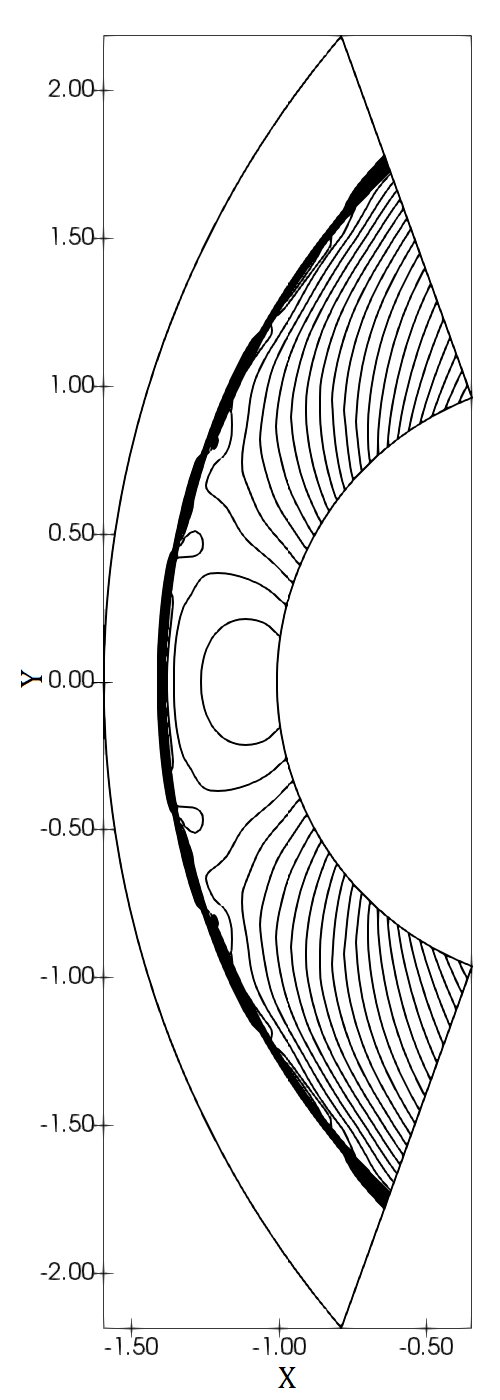} \\
	(a) HLLE  \hspace{12mm} (b)HLLEM  \hspace{12mm} (c) HLLEM-FP1D $r=\frac{1}{2}$ \hspace{2mm} (d) HLLEM-FP1D  $r=\frac{1}{3}$
	\caption{Density contours for Mach 20 flow over a blunt body computed by the HLLE, HLLEM and HLLEM-FP1D schemes. The results are shown after 100,000 iterations}
	\label{carbuncle-hllem1d}
	\end{center}
\end{figure}

\begin{figure} [H]
	\begin{center}
	\includegraphics[width=300pt]{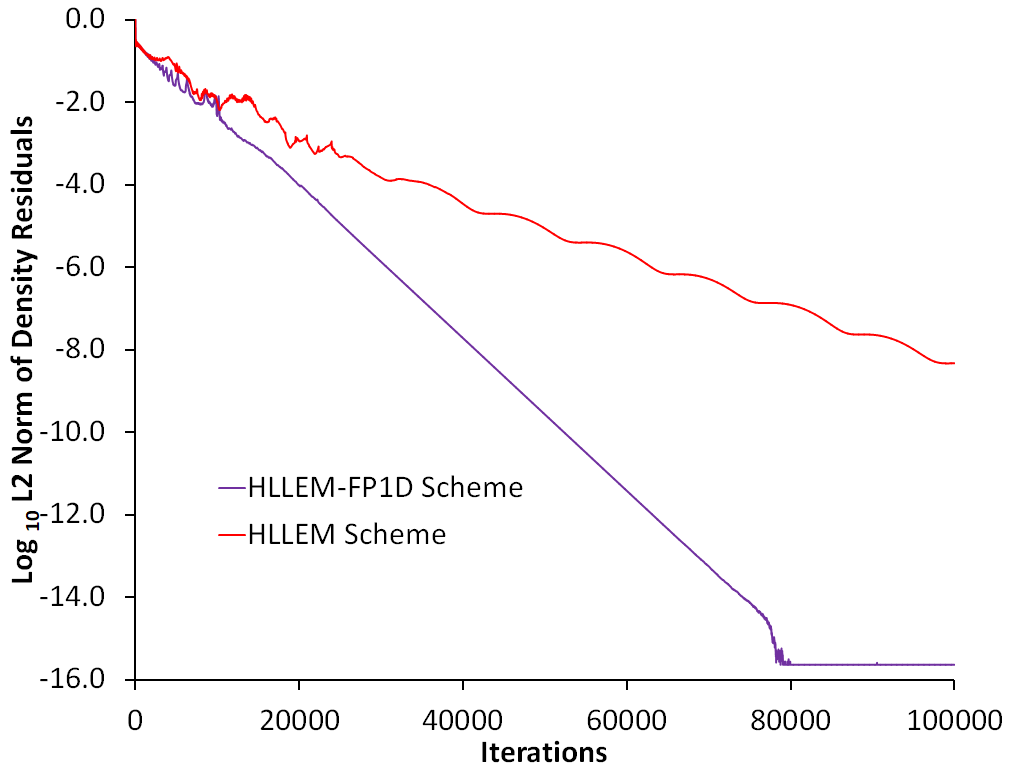} 
	\caption{History of L2 norm of density residuals for $M_\infty=20$ flow past a blunt body computed by the HLLEM and HLLEM-FP1D schemes.} 
	\label{convergence-carbuncle-hllem-fp1d}
	\end{center}
\end{figure}

Hypersonic flow over a blunt body is a typical problem to assess the performance of a scheme for the carbuncle instability. A free-stream Mach number of 20 is considered in the computation. The grid for the blunt body had a size of $40\times320$ cells. The domain is initialized with values of $\rho=1.4, p=1, u=20$ and $v=0$. The inlet boundary condition is set to free-stream values. The solid wall boundary condition is applied to the blunt body. The computations are carried out for 100,000 iterations. The contour plot for density for the original and the modified HLLEM schemes are shown in Fig. \ref{carbuncle-hllem1d} and a total of 27 density contours from 2.0 to 8.7 are drawn. It can be seen from the figure that the original HLLEM scheme exhibits a severe carbuncle phenomenon, while the proposed HLLEM-FP1D scheme is free from the carbuncle problem of the original scheme. The convergence history plots of the original and proposed HLLEM schemes are shown in Fig. \ref{convergence-carbuncle-hllem-fp1d}. It can be seen from the figure that the original HLLEM scheme converges to about $10^{-8}$ while the proposed HLLEM-FP1D scheme converges to machine accuracy within 100,000 iterations.
\subsection{Supersonic corner problem}
\begin{figure}[H]
	\begin{center}
	\includegraphics[width=220pt]{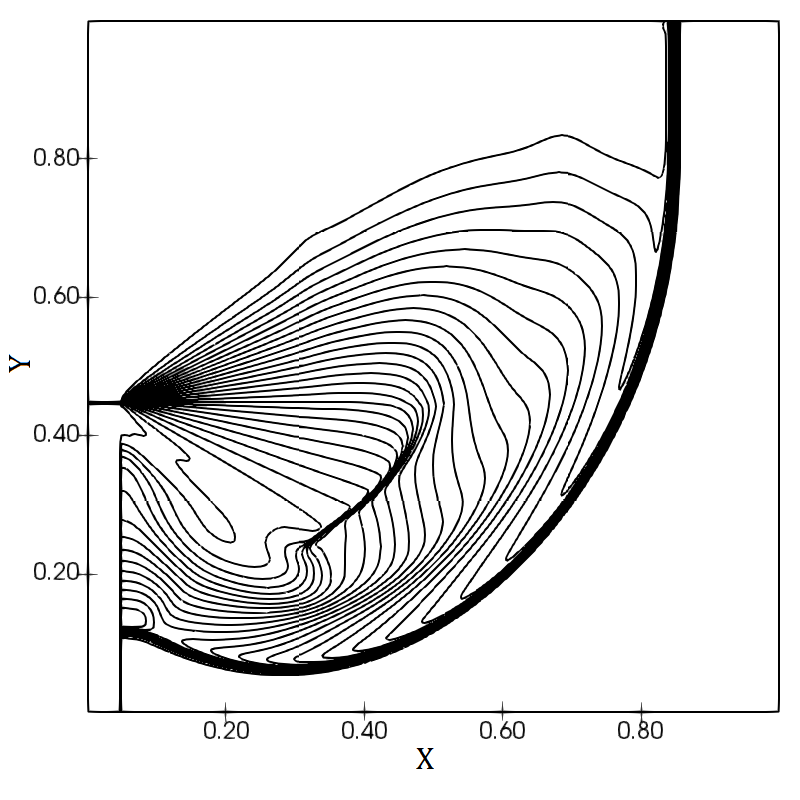} 	\includegraphics[width=220pt]{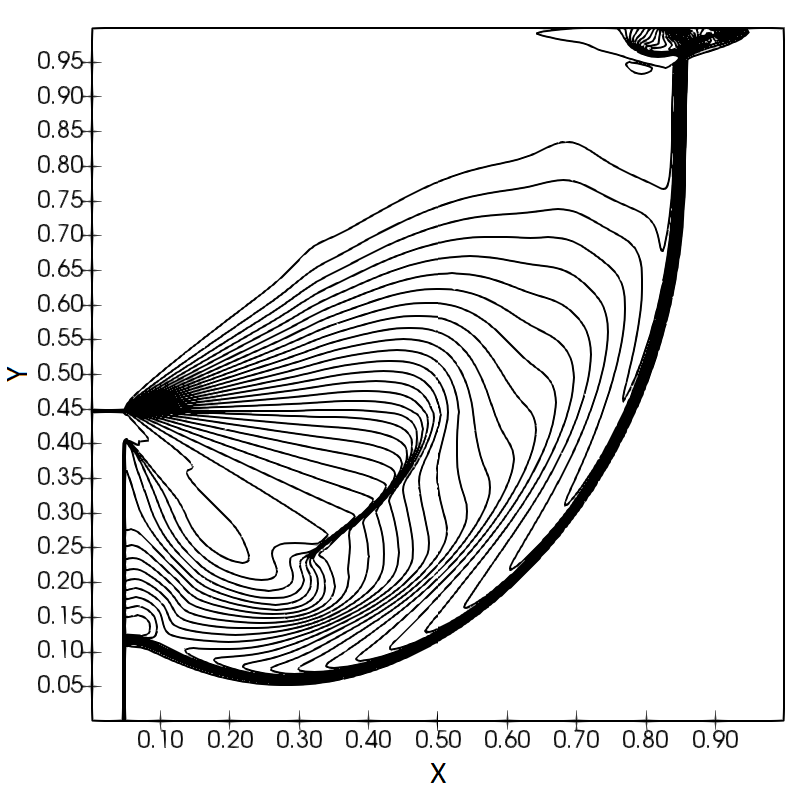} \\
	(a) HLLE Scheme \hspace{3cm} (b)  HLLEM Scheme  \\
	\includegraphics[width=220pt]{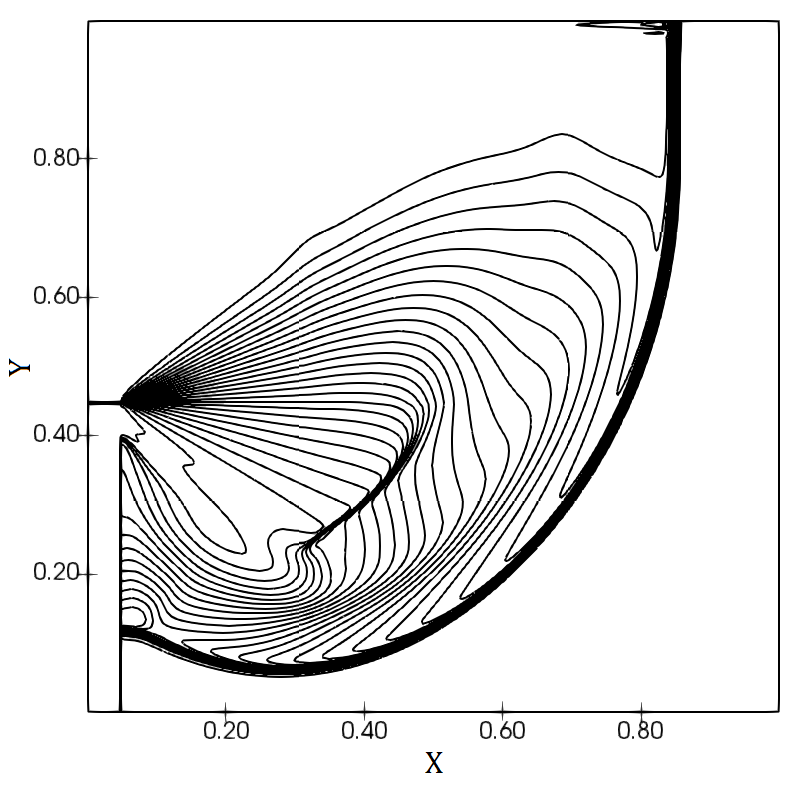} \includegraphics[width=220pt]{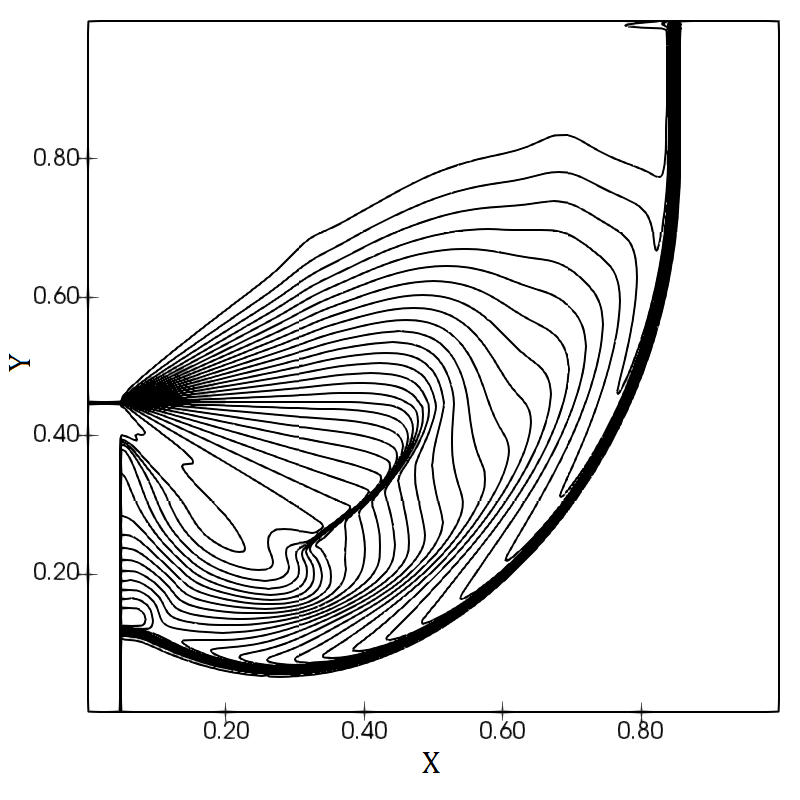} \\ 
	(c) HLLEM-FP1D Scheme $r=\frac{1}{2}$\hspace{3cm} (b) HLLEM-FP1D Scheme $r=\frac{1}{3}$ 
	\caption{Density contours for $M_{\infty}=5.09$ normal shock diffraction around a $90^0$ corner computed by the HLLE, HLLEM and HLLEM-FP1D schemes. The results are shown for time $t=0.1561$.}
	\label{scr-hllem1d}
	\end{center}
\end{figure}

The problem consists of a sudden expansion of Mach 5.09 normal shock around a 90-degree corner. The domain is a square of one unit and is divided into $400\times400$ cells. The corner is located at x=0.05 and y=0.45. The initial normal shock is located at $x=0.05$. The domain to the right of shock is assigned initially with pre-shock conditions of $\rho{}=1.4$, $p=1.0$, $u=0.0$, $v=0.0$. The domain to the left of the shock is assigned post-shock conditions. The inlet boundary is supersonic, the outlet boundary has zero gradients and the bottom boundary behind the corner uses extrapolated values. Reflective wall boundary conditions are imposed on the corner. A CFL value of 0.8 is used and the results are generated for time $t=0.1561$ units. A total of 30 density contours ranging from 0 to 7.1 is shown in Fig. \ref{scr-hllem1d} for the HLLE, HLLEM and the proposed HLLEM schemes. It can be seen from the figure that the classical HLLE scheme is free from numerical shock instability while a severe shock instability is present in the original HLLEM scheme. The proposed HLLEM-FP1D scheme is free from the numerical shock instability problem of the original scheme. 

\section{Conclusions}
This paper utilizes the direct Lyapunov method to assess the stability of HLL-family schemes. Furthermore, we provide an analysis that incorporates these findings to formulate a shock-stable HLLEM scheme. the key conclusions drawn from this study are as follows
\begin{enumerate}
\item A novel approach of non-linear stability analysis based on direct Lyapunov method for numerical shock instability problem in approximate Riemann solvers is proposed here. This is a general approach that can be applied to any non-linear Riemann solver. Its application to several HLL-type schemes is demonstrated here.
\item The HLLE scheme is globally, asymptotically stable whereas the HLLEM scheme is not asymptotically stable. The pressure perturbation feeding the density and transverse momentum perturbations is identified as the cause of instability in the Roe, HLLC and HLLEM schemes. The present methodology offers a more detailed understanding of the origins of numerical shock instability within the Roe, HLLEM and HLLC schemes, surpassing the insights provided by Quirk's linear perturbation method. Stability evaluation using the direct Lyapunov method indicate that the Roe, HLLC and HLLEM schemes are simple stable in the absence of pressure perturbation. The present stability evaluation also indicate that the magnitude of the instabilities in these schemes depends upon the magnitude of the pressure perturbations. The magnitude of the instabilities can therefore be reduced by reducing the magnitude of the pressure perturbations. Available literature on the approximate Riemann solvers suggests that the magnitude of the pressure perturbations in the plane transverse to the normal shock can be reduced by an appropriate low Mach fix. Numerical experiments on the HLLEM scheme validates the link between the magnitude of numerical shock instabilities and magnitude of pressure perturbations discovered by the present stability analysis.
\item The contact-capturing HLLCM/HLLEC schemes and the shear-resolving HLLS/HLLES schemes are not asymptotically stable. The pressure perturbation feeding the density and transverse momentum perturbations is identified as the cause of instability in the contact-capturing HLLCM/HLLEC schemes. The density perturbation feeding the transverse momentum perturbation is identified as the cause of instability in the shear-resolving HLLS/HLLES scheme.
\item The understanding of instability behaviour from direct Lyapunov method based analysis can be used to construct shock stable version of any Riemann solver, which has been demonstrated for HLLEM scheme as an example. The shock-stable HLLEM scheme is obtained by introducing damping of the density and transverse momentum perturbations in presence of pressure perturbations and by reducing the feeding of the pressure perturbations to the density and transverse momentum perturbations. The magnitude of the pressure perturbations in the plane transverse to the grid-aligned shock is reduced by the low Mach fix of Dellacherie for Godunov-type schemes. Numerical results on several test cases show that the proposed HLLEM-FP1D scheme is free from the numerical shock instability problems of the original HLLEM scheme and is also capable of resolving the contact discontinuities accurately at the same time. Similar approach may be attempted for developing shock-stable versions of other complete Riemann solvers like Roe and HLLC.
\end{enumerate}
\end{doublespace}

\end{document}